\documentclass[11pt]{article}
\usepackage{amssymb}
\usepackage{amsmath}

\newtheorem{theo}{Theorem}[section]
\newtheorem{prop}[theo]{Proposition}
\newtheorem{lemm}[theo]{Lemma}
\newtheorem{coro}[theo]{Corollary}
\newtheorem{rema}[theo]{Remark}

\newtheorem{conj}[theo]{Conjecture}

\newcommand{\cqfd}
{%
\mbox{}%
\nolinebreak%
\hfill%
\rule{2mm}{2mm}%
\medbreak%
\par%
}
\newfont{\gothic}{eufb10}
\date{\empty}
\begin{document}
\title{On the Chow ring of certain algebraic \\hyper-K\"ahler manifolds}
\author{Claire Voisin\\ Institut de math\'{e}matiques de Jussieu, CNRS,UMR
7586} \maketitle \setcounter{section}{-1}
\begin{flushright} {\it Pour  Fedya Bogomolov, en l'honneur de ses 60 ans}
\end{flushright}
\section{Introduction}
This paper proposes and studies  a generalization of a conjecture made by Beauville in \cite{beauBB}.
Recall that Beauville and the author proved the following result in \cite{BeVo}:
\begin{theo} \label{bv1}Let $S$ be an algebraic $K3$ surface. Then there exists a degree $1$
$0$-cycle $o$ on $S$ satisfying the property that for any line bundle $L$ on $S$,
one has
$$c_1(L)^2=[c_1(L)]^2o\,\,{\rm in}\,\, CH_0(S).$$
Furthermore, we have $c_2(T_S)=24 o$.
\end{theo}
(In this paper,  Chern classes will be  Chern classes in the Chow
ring tensored by $\mathbb{Q}$, and we will denote by $[c_i]$ the
corresponding rational cohomology classes.)

This result can be rephrased by saying that any polynomial relation
$$P([c_1(L_i)])=0\,\,{\rm in}\,\,H^*(S,\mathbb{Q}),\,\, L_i\in Pic\,S,$$
already holds in $CH(S)$.

In \cite{beauBB}, Beauville conjectured that a similar result holds for algebraic hyper-K\"ahler
varieties:
\begin{conj}\label{conj1} (Beauville) Let $Y$ be an algebraic hyper-K\"ahler
variety. Then any polynomial cohomological relation
$$P([c_1(L_i)])=0\,\,{\rm in}\,\,H^*(Y,\mathbb{Q}),\,L_i\in Pic\,Y$$
already holds at the level of Chow groups :
$$P(c_1(L_i))=0\,\,{\rm in}\,\,CH(Y).$$
\end{conj}
He proved in \cite{beauBB} this conjecture in the case of the second and third punctual Hilbert scheme
of an algebraic $K3$ surface.

In this paper, we observe that the results of \cite{BeVo} can lead
to a more general conjecture concerning the Chow ring
of an algebraic hyper-K\"ahler
variety. Namely, the full statement of Theorem \ref{bv1}
 can  be interpreted by saying that any polynomial relation
between $[c_2(T_S)], [c_1(L_i)]$ in $H^*(S,\mathbb{Q})$, already holds
between
$c_2(T_S),\,c_1(L_i)$ in $CH(S)$.
The purpose of this paper is to study the following conjecture:
\begin{conj}\label{conj2}Let $Y$ be an algebraic hyper-K\"ahler
variety. Then any polynomial cohomological relation
$$P([c_1(L_j)], [c_i(T_Y))])=0\,\,{\rm in}\,\,H^{2k}(Y,\mathbb{Q}),\,L_j\in Pic\,Y$$
already holds at the level of Chow groups :
$$P(c_1(L_j), c_i(T_Y))=0\,\,{\rm in}\,\,CH^k(Y).$$
\end{conj}
We shall prove the following results:
\begin{theo}\label{maintheo} 1)  Conjecture \ref{conj2} holds for $Y=S^{[n]}$, for
$n\leq 2b_{2}(S)_{tr}+4$ and any $k$, where $S^{[n]}$ is the Hilbert scheme of length $n$
subschemes of an algebraic $K3$ surface $S$.

2) Conjecture \ref{conj2} is true  for  any $k$ when $Y$ is the
Fano variety of lines of a cubic fourfold.
\end{theo}
In 1), $b_{2}(S)_{tr}=b_2(S)-\rho$
is the rank of the transcendental lattice of $S$.

Concerning point 2), recall from \cite{beaudonagi} that the  variety of lines $F$ of a
cubic fourfold $X$ is a deformation of $S^{[2]}$, for $S$ an
algebraic $K3$ surface, but that for general $X$, it has
$Pic\,F=\mathbb{Z}$ and thus it is not a Hilbert scheme. Even when
$\rho(F)\geq2$ it is not necessarily the case that $F$ is a
$S^{[2]}$. In \cite{beauBB}, Beauville asked whether his conjecture
\ref{conj1} holds true for the variety of lines  of a cubic
fourfold.

Finally, we also prove the following.
\begin{theo}\label{smalldim} Conjecture \ref{conj2} holds  for $Y=S^{[n]}$, and $k=2n-2,\,2n-1,\,2n$, for
any $S$ as above and any $n$.
\end{theo}

The cohomology ring of the Hilbert scheme of a $K3$ surface has been computed in \cite{fango}, \cite{LS}.
 For the subring generated by $H^2$, on can use the result of Verbitsky \cite{Ver}, \cite{Bo}.
 The question of understanding more precisely the Chow ring is rather delicate and we are dealing here
 only with a small part of it.

We prove in section 1 part 1) of Theorem \ref{maintheo}
 and Theorem \ref{smalldim} . The proof involves particular cases of the following statement
:
\begin{conj}\label{conj3} Let $S$ be an algebraic $K3$ surface. For any integer $m$, let
$P\in CH(S^m)$ be a polynomial expression in
$$pr_i^*c_1(L_s),\,L_s\in Pic\,S,\,pr_j^*o,\,pr_{kl}^*\Delta_S.$$
Then if $[P]=0$, we have $P=0$.
\end{conj}
We also prove that  Conjecture \ref{conj3} for $S$ and any $m'\leq
m$ implies Conjecture \ref{conj2} for $Y=S^{[m]}$.

 In section 2 we deal with the case of  the variety of
lines of the cubic fourfold (Theorem \ref{maintheo}, 2)).

It is a pleasure to dedicate this paper to Fedya Bogomolov, who greatly  contributed
in the papers \cite{bo1}, \cite{bo2}, \cite{Bo} to the study of hyper-K\"ahler manifolds.

\section{Case of the Hilbert scheme of a $K3$ surface}
Let $S$ be an algebraic  $K3$ surface, and $S^{[n]}$ be the Hilbert
scheme of length $n$ subschemes of $S$. For any line bundle $L$ on
$S$, there is an induced line bundle, which we still denote by $L$
on $S^{[n]}$, which is the pull-back via the Hilbert-Chow morphism
of the line bundle on $S^{(n)}$ corresponding to the
$\mathfrak{S}_n$-invariant line bundle $L\boxtimes\ldots\boxtimes L$
on $S^n$.

There are furthermore two natural vector bundles on $S^{[n]}$,
namely $\mathcal{O}_{[n]}$, which is defined as
$R^0p_*\mathcal{O}_{\Sigma_n}$, where
$$\Sigma_n\subset S^{[n]}\times S,\,{p=pr_1}:\Sigma_n{\rightarrow}S^{[n]}$$
is the incidence scheme, and the tangent bundle $T_n$. It is not
clear that the Chern classes of $\mathcal{O}_{[n]}$ can be expressed
as polynomials in $c_1(\mathcal{O}_{[n]})$ and the Chern classes of
$T_n$. The following result may thus be stronger than Theorem
\ref{maintheo}, 1):
\begin{theo} \label{theo2}Let $n\leq 2b_2(S)_{tr}+4$, and let
$P\in CH(S^{[n]})$ be any polynomial expression  in the variables
$$c_1(L),\, L\in Pic\,S\subset Pic\,S^{[n]}, \, c_i(\mathcal{O}_{[n]}),\,c_j(T_n)\in CH(S^{[n]}).$$ Then if
$P$ is cohomologous to $0$, we have $P=0$ in $CH(S^{[n]})$.
\end{theo}

This implies Theorem \ref{maintheo} for the $n$-th Hilbert scheme
of $K3$ surface $S$ with $n\leq 2b_2(S)_{tr}+4$, because we have
$c_1(\mathcal{O}_{[n]})=-\delta$, where $2\delta\equiv E$ is the
class of the  exceptional divisor of the resolution
$S^{[n]}\rightarrow S^{(n)}$, and it is well-known that
$Pic\,S^{[n]}$ is generated by $Pic\,S$ and $\delta$.

To start the proof of this theorem, we establish first the following
Proposition \ref{prop1}, which gives particular cases of Conjecture
\ref{conj3}.
 Let $o\in CH^2(S)$ be the cycle introduced in the introduction. Let $m$ be an integer.
\begin{prop}\label{prop1} Let $P\in CH(S^m)$ be a polynomial expression in the variables
$$pr_i^*(\frac{1}{24}c_2(T))= pr_i^*o,\,pr_j^* c_1(L_s),\,L_s\in Pic\,S,\,pr_{kl}^*\Delta_S,\,k\not=l,$$
where $\Delta_S\subset S\times S$ is the diagonal. Assume that one
of the following assumptions is satisfied: \begin{enumerate}
\item $m\leq 2b_{2}(S)_{tr}+1$.
\item  $ P$ is invariant under the action of the
symmetric group $ \mathfrak{S}_{m-2}$ acting on the $m-2$ first
indices.
\end{enumerate}
Then if $P$ is cohomologous to
$0$, it is equal to $0$ in $CH(S^m)$.
\end{prop}

Using the results of \cite{BeVo}, this proposition is a consequence
of the following lemma:
\begin{lemm}\label{lemm1}
The polynomial relations $[P]=0$ in the cohomology ring $H^*(S^m)$,
satisfying one of the above assumptions  on $m,\,P$, are all
generated (as elements of  the ring of all polynomial expressions in
the variables above) by the following polynomial relations, the list
of which will be denoted by (*) :
\begin{enumerate}
\item \label{item1} $[pr_i^*(c_1(L))\cdot pr_i^*o]=0,\,L\in Pic\,S$,  $[pr_i^*(o)\cdot pr_i^*(o)]=0$.
 \item $[pr_i^*(c_1(L)^2-[c_1(L)]^2o)]=0,\,L\in Pic\,S$.
 \item \label{item3} $[pr_{ij}^*(\Delta_S.p_1^*o-(o,o))]=0$, where $p_1$ here is the first projection
of $S\times S$ to $S$, and $(o, o)=p_1^*o\cdot p_2^*o$.
\item $[pr_{ij}^*(\Delta_S.p_1^*c_1(L)-c_1(L)\times o-o\times c_1(L))]=0,\,L\in Pic\,S$,
where $p_1$ here is the first projection of $S\times S$ to $S$, and
$c_1(L)\times o=p_1^*c_1(L)\cdot p_2^*o$.
\item \label{5} $[pr_{ijk}^*(\Delta_3-p_{12}^*\Delta_S\cdot p_3^* o-p_1^*o\cdot p_{23}^*\Delta_S-
p_{13}^*\Delta_S\cdot p_2^*o+
p_{12}^*(o,o)+p_{23}^*(o,o)+p_{13}^*(o,o))]=0$.
\item $[pr_{ij}^*\Delta_S]^2=24pr_{ij}^*(o,o)=24pr_i^*o\cdot pr_j^*o$.
\end{enumerate}
\end{lemm}
In \ref{5}, $\Delta_3$ is the small diagonal of $S^3$
and the $p_i,\,p_{ij}$ are the various projections from $S^3$ to
$S$, $S\times S$ respectively. Note that $\Delta_3$ can be expressed
as $p_{12}^*\Delta_S\cdot p_{23}^*\Delta_S$. Furthermore we have
$$pr_{ij}^*\circ p_1^*=pr_i^*,\,pr_{ijk}^*\circ p_{12}^*=pr_{ij}^*,\,pr_{ijk}^*\circ p_i^*=pr_i^*.$$
  Thus all the relations in (*)
are  actually polynomial expressions in the variables
$$[pr_i^*o],\,[pr_j^* c_1(L)],\,L\in
Pic\,S,\,[pr_{kl}^*\Delta_S],\,k\not=l.$$

Assuming this Lemma, we conclude that for $m \leq 2b_{2}(S)_{tr}+1$, all  polynomial
relations $[P]=0$ in the variables $pr_i^*o,\,pr_j^*c_1( L),\,L\in
Pic\,S,\,pr_{kl}^*\Delta_S,\,k\not=l$ which hold in $H^*(S^m)$ also
hold in $CH(S^m)$, because we know from \cite{BeVo} that the
relations listed in (*) hold in $CH(S^m)$. In fact, (apart from the
relations \ref{item1} and \ref{item3} which obviously hold in
$CH(S^m)$), these relations are pulled-back, via the maps $pr_{i}$,
resp. $pr_{ij}$, resp. $pr_{ijk}$, from relations in $CH(S)$, resp.
$CH(S^2)$, resp. $CH(S^3)$, which are established in \cite{BeVo}.

Similarly, for any $m$, the same conclusion holds for polynomial
relations invariant under $\frak{S}_{m-2}$.

This concludes the proof of Proposition \ref{prop1}.\cqfd

{\bf Proof of Lemma \ref{lemm1}.} Let ${\mathcal B}$ be a basis of
$Pic\,S$. It is clear that modulo the relations generated by (*),
any  polynomial in the variables
\begin{eqnarray}\label{varmentonlundi}[pr_{ij}^*\Delta_S],\,[pr_k^*c_1(L)],\,L\in {\mathcal B},\,
[pr_l^*o],
\end{eqnarray}
 can be written as a
combination of monomials having the property that an index $i\in
\{1,\ldots,n\}$ appears only once. Indeed, these relations express
any product with a repeated index as a combination of monomials with
no repeated index. Furthermore, if we start from a polynomial which
is invariant under the action of $\frak{S}_{m-2}$, as the set of
relations (*) is stable under this action, it is clear that
replacing systematically each repeated index by the corresponding
combination with no repeated index using (*), we will end with a
 polynomial expression invariant under the action of $\frak{S}_{m-2}$.

We claim now that  if $m\leq 2b_{2}(S)_{tr}+1$, no non zero combination of monomials
with no repeated index vanishes in $H^*(S^m)$. Furthermore, for any
$m$,  no non zero  combination of monomials with no repeated index
which is invariant under $\frak{S}_{m-2}$  vanishes in $H^*(S^m)$.

To prove the claim, consider   the transcendental part of $H^2(S,{\mathbb Q})$,
$$H^2(S,{\mathbb Q})_{tr}:=NS(S)^{\perp}. $$
We have the direct sum decomposition
\begin{eqnarray}\label{direct}
H^*(S,{\mathbb Q})=H^2(S,{\mathbb Q})_{tr}\oplus H^*(S,{\mathbb Q})_{alg},
\end{eqnarray}
where $H^*(S,{\mathbb Q})_{alg}$ is generated by $$H^0(S,{\mathbb
Q}),\,NS(S)_{\mathbb Q},\,H^4(S,{\mathbb Q}).$$
 The
decomposition (\ref{direct}) induces for any $m$ a direct sum
decomposition of
$$H^*(S^m,{\mathbb Q})=H^*(S,{\mathbb Q})^{\otimes m}=pr_1^*H^*(S,{\mathbb Q})\otimes\ldots
\otimes pr_m^*H^*(S,{\mathbb Q}).$$ Let $[\Delta_S]_{tr}$ be the
projection of $[\Delta_S]$ in the direct summand
$$H^2(S,{\mathbb Q})_{tr}\otimes H^2(S,{\mathbb Q})_{tr}$$
of $H^4(S\times S, {\mathbb Q})$. Then we have
$[\Delta_S]_{tr}\not=0$ and $$
[\Delta_S]=[\Delta_S]_{tr}+p_1^*[o]+p_2^*[o]+\sum_{i,j\in {\mathcal
B}}\alpha_{ij}p_1^*[c_1(L_i)]\cdot p_2^*[c_1(L_j)].$$ It is then
clear that it suffices to prove the claim with $pr_{ij}^*[\Delta_S]$
replaced by $pr_{ij}^*[\Delta_S]_{tr}$ in the set of variables
(\ref{varmentonlundi}).

 Let $M$ be a
monomial of the form above, and let $I_M\subset \{1,\ldots,m\}$ be
the set of indices $i$ appearing in $M$ via a diagonal, i.e. for
some $l\not=i$, the variable $pr_{il}^*[\Delta_S]_{tr}$ appears in
$M$. Then $I_M$ is also the unique index  set for which the
projection of $M$ in
$$\bigotimes_{i\not \in I_M}H^*(S,{\mathbb Q})_{alg}\otimes \bigotimes_{i\in I_M}
H^2(S,{\mathbb Q})_{tr}$$ is non zero. Hence it follows that a
relation
$$\sum_M\alpha_M M=0$$
implies for each  fixed $I\subset \{1,\ldots,m\}$ (of even
cardinality), by projection onto $$\bigotimes_{i\not\in
I}pr_i^*H^*(S,{\mathbb Q})_{alg}\otimes \bigotimes_{i\in I}pr_i^*
H^2(S,{\mathbb Q})_{tr},$$ a relation of the form
\begin{eqnarray}\label{clo}\sum_{M,\,I_M=I}\alpha_M M=0.
\end{eqnarray}

Now note that we can further decompose each term
$\bigotimes_{i\not\in I}pr_i^*H^*(S,{\mathbb Q})_{alg}$ using the
basis of $H^*(S,\mathbb{Q})_{alg}$ given by $H^4(S,{\mathbb
Q}),\,H^0(S,{\mathbb Q})$ and the basis ${\mathcal B}$. Then the
relation (\ref{clo}) decomposes into a sum of relations of the form
\begin{eqnarray}\label{relinutile}(\otimes_{j\not\in I}pr_j^*\alpha_j)\otimes
(\sum_{I_{M'}=I}\alpha_{M'}' M'), \end{eqnarray} for any given $I$
and given set $(\alpha_j)_{j\not\in I}$ of chosen elements in the
basis
$$H^4(S,{\mathbb Q}),\,H^0(S,{\mathbb Q}),\,{\mathcal B}$$ of
$H^*(S,{\mathbb Q})_{alg}$. Here each $M'$ is a monomial in the
$pr_{ij}^*[\Delta_S]_{tr}$, and $I_{M'}=I$ means that only indices
$i,\,j\in I$ appear in the monomial $M'$, and each index $l\in I$
appears exactly once.

Of course, (\ref{relinutile})  is equivalent to the relation
\begin{eqnarray}\sum_{I_{M'}=I}\alpha_{M'}' M'=0,
\label{relmermenton}
\end{eqnarray}
which has to hold in $H^{2s}(S^{m},{\mathbb Q})$ or equivalently in
$H^{2s}(S^{I},{\mathbb Q})$. (Here $2s$ is the cardinality of $I$,
and $S^I$ is the product of the copies of $S$ indexed by $I$. Thus
clearly (\ref{relmermenton}) is pulled-back via the projection
$S^m\rightarrow S^I$ from the corresponding relation in $S^I$.)

Now observe that if we started with a  polynomial relation invariant
under the action of $\frak{S}_{m-2}$, each relation we get in
(\ref{relmermenton}) is invariant under the symmetric group
$\mathfrak{S}_{I}'$ permuting the elements of $I$ which are $\leq
m-2$.

In conclusion, we are reduced to prove that for each $I\subset \{1,\ldots,m\}$ of cardinality
$2s$, and thus satisfying
$2s\leq m\leq 2b_{2}(S)_{tr}+1$, there are no  relations in $H^{2s}(S_I)$ between
monomials of the form
\begin{eqnarray}\label{monomial}\prod pr_{ij}^*[\Delta_S]_{tr},\end{eqnarray}
where each index $i,\,j\in  I$ appears exactly once. Furthermore,
for any $m$,  there are no  $\frak{S}_I'$-invariant relations in
$H^{2s}(S_I)$ between monomials of the form above.

For the  statement concerning $\frak{S}_I'$-invariant relations,
 this is
obvious, as the symmetric group $\mathfrak{S}_I'$ acts with at most
two distinct orbits on the set of monomials of degree $s$ in the variables
$pr_{ij}^*[\Delta_S]_{tr},\,i,\,j\in I,$ with no repeated indices,
namely, in the case where $m-1,\,m\in I$, those monomials containing
$p_{m-1,m}^*[\Delta_S]_{tr}$ and those not containing it. Assume
that $m-1,\,m\in I$, as otherwise the action is transitive and the
result is still simpler.

  Then the only possible non zero  $\frak{S}_I'$-invariant relation
would be of the form: \begin{eqnarray}
\label{finalmermenton}\beta\sum_{M\in{\mathcal E}}M
=\alpha\sum_{M\in{\mathcal F}}M,
\end{eqnarray} where ${\mathcal E}$ is the set of  monomials not containing
$pr_{m-1,m}^*[\Delta_S]_{tr}$, and ${\mathcal F}$ is the set of
monomials  containing $pr_{m-1,m}^*[\Delta_S]_{tr}$. We identify $I$
with $\{1,\ldots,2s\}$ in such a way that $m-1$ is identified with
$2s-1$ and $m$ with $2s$. Then this gives an identification of $S^I$
with $S^s\times S^s$ and we can consider each $M$ as above as a
self-correspondence of $S^s$. Let $p_1,\, p_2$ be the two
projections of $S^{2s}$ to $S^s$. Then each monomial $M $ as above
induces a map
$$\gamma_M:H^{2s,0}(S^s)\rightarrow H^{2s,0}(S^s),$$
given by
\begin{eqnarray}\label{10mai}\gamma_M(\eta)=p_{1*}(M\cdot p_2^* \eta).
\end{eqnarray}
Now we have $\gamma_M=0$ when the monomial  $M$
has the property that for two indices $i,\,j\leq s$,
$pr_{ij}^*[\Delta_{S}]_{tr}$ divides $M$, because $\omega^2=0$ on $S$.
In
particular, we have $\gamma_M=0$ for $M\in {\mathcal F}$. Thus
(\ref{finalmermenton}) gives
\begin{eqnarray}\label{25oct}\beta\sum_{M\in {\mathcal
E}}\gamma_M=0.\end{eqnarray}

 On the other hand, when   the
monomial $M$
has the property that for any indices $i\not=j\leq s$,
$pr_{ij}^*[\Delta_{S}]_{tr}$ does not divide  $M$, $M$ is of the form
$$\Pi_{i\leq s}pr_{i, s+\sigma(i)}^*[\Delta_S]_{tr}$$
for some permutation $\sigma$ of $\{1,\ldots,s\}$. In that case, we
have \begin{eqnarray}
\label{formval16}\gamma_M(pr_1^*\omega\cdot\ldots\cdot pr_s^*\omega)=
pr_{\sigma(1)}^*\omega\cdot\ldots\cdot
pr_{\sigma(s)}^*\omega=p_1^*\omega\cdot\ldots\cdot p_s^*\omega.
\end{eqnarray}
  Thus we get by (\ref{25oct}) and (\ref{formval16})
$$\beta\sum_{M\in {\mathcal
E}}\gamma_M=\beta s!Id_{H^{2s,0}(S^s)}=0,$$ which implies that
$\beta=0$, and that if $\alpha\not=0$, the relation reduces to $\sum_{M\in {\mathcal
F}}M=0$. But elements of ${\mathcal F}$ are of the form
$$p_{m-1,m}^*[\Delta_S]_{tr}\cdot M',$$
where $M'\in {\mathcal F}'$ are the monomials with no repeated
indices in the $pr_{ij}^*[\Delta_S]_{tr},\,i,\,j<m-1,\,i,\,j\in I$.
The relation $\sum_{M\in {\mathcal F}}M=0$ thus provides $\sum_{M'\in
{\mathcal F}'}M'=0$, which has just been proved to be impossible.

In the case where $2s\leq  m\leq 2b_{2}(S)_{tr}+1$, we have  $s\leq b_{2}(S)_{tr}$.
The intersection form on $H^2(S,{\mathbb C})_{tr}$ is non degenerate. Thus we can
choose an orthonormal  basis $\alpha_i,\,1\leq i\leq b_{2}(S)_{tr},$ of $H^2(S,{\mathbb C})_{tr}$.
 Let
$$\eta:=pr_1^*\alpha_1\cdot\ldots\cdot pr_s^*\alpha_s.$$
For each monomial $M$, consider now
$\gamma_M: H^{2s}(S^s,\mathbb{C})\rightarrow H^{2s}(S^s,C)$, defined as in (\ref{10mai}).
As we have $<\alpha_i,\alpha_j>=0$ for $i\not=j$, we find that
if there are two indices $i,\,j> s$ such that $p_{ij}^*[\Delta_S]_{tr}$ appears in $M$, we have
$\gamma_M(\eta)=0$, and that the remaining  $M$'s are in one-to-one correspondence
with  permutations $\sigma$ of
$\{1,\ldots,s\}$. Then, as before, we have for such $M$:
$$\gamma^*(\eta)=p_{\sigma(1)}^*\alpha_1\cdot\ldots\cdot p_{\sigma(s)}^*\alpha_s=
p_1^*\alpha_{\sigma^{-1}(1)}\cdot\ldots\cdot p_{s}^*\alpha_{\sigma^{-1}(s)}.$$

As the tensors $
pr_1^*\alpha_{\sigma^{-1}(1)}\cdot\ldots\cdot pr_s^*\alpha_{\sigma^{-1}(s)}$, $\sigma\in\mathfrak{S}_s$
are linearly independent in $H^2(S,{\mathbb C})_{tr}$, we conclude from these two facts that a relation
$\sum_M\alpha_MM=0$ implies $\alpha_M=0$ for all those $M$ such that for no
indices $i,\,j> s$, $pr_{ij}^*[\Delta_S]_{tr}$ appears in $M$.

To show that the other coefficients $\alpha_M$ must be also $0$, we introduce maps similar
to the $\gamma_M$, defined by choosing any subset $I_1=\{i_1,\ldots,i_s\}$,
$ i_1<\ldots<i_s$ of $I$. Denoting by $I_2$ the complementary set, we define
$\gamma_M^{I_1}$ as $p_{I_1*}\circ (M\cdot p_{I_2}^*)$, where
$p_{I_1}$ (resp. $p_{I_2}$)  is the  projection from $S_I\cong S^{2s}$ to $S^s$ determined
by the (ordered) set $I_1$ (resp. $I_2$).
For any $M$, there is a choice of $I_1$ such that
 for no
indices $i,\,j\in I_2$, $pr_{ij}^*[\Delta_S]_{tr}$ appears in $M$, and then we conclude as before that
$\alpha_M$ must also be $0$.

Thus Lemma \ref{lemm1} and also Proposition \ref{prop1} are proved.

\cqfd
We come now to the geometry of $S^{[n]}$.
Let us introduce   the following notation: let
$$\mu=\{\mu_1,\ldots,\mu_m\},\,m=m(\mu),\,\sum_i\mid \mu_i\mid =n$$
be a partition of $\{1,\ldots,n\}$. Such a partition determines a
partial diagonal
$$S_\mu\cong S^m\subset S^n,$$
defined by the conditions
$$x=(x_1,\ldots, x_n)\in S_\mu
\Leftrightarrow x_i=x_j\,\,{\rm if}\,\, i,\,j\in \mu_l,\,\,{\rm for\,\,some}\,\,l.$$

Consider the quotient map
$$q_\mu:S^m\cong S_\mu\rightarrow S^{(n)},$$
and denote by $E_\mu$ the following fibered product:
$$E_\mu:=S_\mu\times_{S^{(n)}} S^{[n]}\subset S^m\times S^{[n]}.$$
We view $E_\mu$ as a correspondence between $S^m$ and $S^{[n]}$ and
we will denote as usual by $E_\mu^*:CH(S^{[n]})\rightarrow CH(S^m)$
the map
$$\alpha\mapsto{pr_1}_*(pr_2^*(\alpha)\cdot E_\mu).$$
Let us denote by $\frak{S}_\mu$ the subgroup
 of $\frak{S}_m$
 permuting only the indices $i,\,j$ for which the cardinalities  of $\mu_i,\,\mu_j$ are equal.
  The group $\frak{S}_\mu$ can be seen as  the quotient of the global stabilizer
 of $S_\mu$ in $S^n$ by its pointwise stabilizer. In this way the action of
$\frak{S}_\mu$ on $S_\mu\cong S^m$ is induced by the action of
$\frak{S}_n$ on $S^n$.

We have the following result:
\begin{prop}\label{prop2} Let $P\in CH(S^{[n]})$ be a polynomial expression in
$c_i(\mathcal{O}_{[n]}),\,c_j(T_n)$. Then for any $\mu$ as above,
$E_\mu^*(P)\in CH(S^m)$ is a polynomial expression in
$pr_s^*o,\,pr_{lk}^*\Delta_S$. Furthermore, $E_\mu^*(P)$ is
invariant under the group $\frak{S}_\mu$.
\end{prop}
Note that the last statement is obvious, since
 $\frak{S}_\mu$  leaves
invariant the correspondence $E_\mu\subset S_\mu\times S^{[n]}$.

We postpone the proof of this proposition and conclude the proofs of the theorems.

\vspace{0.5cm}

{\bf Proof of Theorem \ref{theo2}}.  From the
work of De Cataldo-Migliorini \cite{Cami}, it follows that the map
$$(E_\mu^*)_{\mu\in Part(\{1,\ldots,n\})}:CH(S^{[n]})\rightarrow \Pi_{\mu} CH(S^{m(\mu)})$$
is injective. Let now $P\in CH(S^{[n]})$ be a polynomial expression
in $c_1(L),\,L\in Pic\,S\subset Pic\,S^{[n]},\, \,
c_i(\mathcal{O}_{[n]}),\,c_j(T_n)\in CH(S^{[n]})$. Note first that
for $L\in Pic\,S$, and for each $\mu$, the restriction of $pr_2^*L$
to $E_\mu\subset S_\mu\times S^{[n]}$ is a pull-back $pr_1^*
{L_\mu}_{\mid E_\mu}$, where $L_\mu\in Pic\, S_\mu= Pic\,S^m$ is
equal to $L^{\otimes\mid\mu_1\mid}\boxtimes\ldots\boxtimes
L^{\otimes\mid\mu_m\mid}$. This follows from the fact that $L$ is
the pull-back of a line bundle on $S^{(n)}$. Note that $L_\mu$ is
invariant under $\frak{S}_\mu$.

Thus it follows from Proposition \ref{prop2} and the projection
formula that for each partition $\mu$, $E_\mu^*(P)$ is a polynomial
expression in $pr_i^*c_1(L),\,pr_k^*o,\,pr_{lm}^*\Delta$ which is
invariant under the group $\frak{S}_\mu$.

 Now, if
$P$ is cohomologous to $0$, each $E_\mu^*(P)$ is cohomologous to
$0$. Let us now
  verify that the assumptions of
Proposition \ref{prop1} are satisfied. Recall that we assume
$n\leq 2b_2(S)_{tr}+4$. If $m(\mu)\leq 2b_2(S)+1$,
 Proposition \ref{prop1} applies. Otherwise, $m(\mu)\geq 2b_2(S)+2$ and,  as $n\leq 2b_2(S)_{tr}+4$, it
follows that the partition $\mu$ contains at most two sets of
cardinality $\geq2$. Thus the group $\frak{S}_\mu$  contains in this case a
group conjugate to $\frak{S}_{m(\mu)-2}$. Proposition \ref{prop1}
thus applies, and gives $E_\mu^*(P)=0$ in $CH(S_\mu)$, for all
$\mu$.

It follows that $P=0$ by the result of De Cataldo-Migliorini. This
concludes the proof of Theorem \ref{theo2}. \cqfd
{\bf Proof of Theorem \ref{smalldim}.} Let $P\in CH^k(S^{[n]})$, with $k\geq 2n-2$ be a
polynomial expression
in $c_1(L),\,L\in Pic\,S\subset Pic\,S^{[n]},\, \,
c_i(\mathcal{O}_{[n]}),\,c_j(T_n)\in CH(S^{[n]})$, and assume that $[P]=0$. Notice that because $k\geq 2n-2$,
we have
$E_\mu^*P=0$ if the image of $E_\mu$ in $S^{[n]}$ has codimension $>2$. This is the case once
$m(\mu)<n-2$. On the other hand, if $m(\mu)\geq m-2$, the partition $\mu$ has at most two sets
$\mu_i$ of cardinality $\geq 2$. Hence for $m(\mu)\geq2$, the group $\frak{S}_\mu$  contains a
group conjugate to $\frak{S}_{m(\mu)-2}$. As $[E_\mu^*P]=0$, and
$E_\mu^*P$ is a $\frak{S}_\mu$-invariant polynomial expression in
$pr_i^*c_1(L),\,pr_j^*o,\,pr_{ij}^*\Delta_S$,  Proposition \ref{prop1}
thus applies, and gives $E_\mu^*(P)=0$ in $CH(S_\mu)$ for $m(\mu)\geq n-2$. As we also have
$E_\mu^*(P)=0$ in $CH(S_\mu)$ for $m(\mu)<n-2$, the theorem of De Cataldo-Migliorini shows that $P=0$.
\cqfd

To conclude, let us notice that Proposition \ref{prop2} and the end
of the proof of Theorem \ref{theo2}   prove the following:
\begin{prop} Conjecture \ref{conj3} for $S$, and any $m\leq n$, implies Conjecture \ref{conj2} for
$S^{[n]}$.
\end{prop}

It remains to prove Proposition \ref{prop2}.
 For the proof, we use the formulas proved in
\cite{EGL}, which allow induction on $n$. As in \cite{EGL}, in order
to get the result by induction, we will need to introduce a more
general induction statement, which is the following: For each
integer $l$, we can also consider the correspondence
$$ E_{\mu,l}:=E_\mu\times\Delta_{ S^l}$$
betweeen $ S_\mu\times S^l$ and $S^{[n]}\times S^l$, where $\Delta_{
S^l}$ is the diagonal of $S^l$. On $S^{[n]}\times S^l$, we have the
natural classes $pr_{0i}^*c_s(\mathcal{I}_n)$, where
$\mathcal{I}_{n}$ is the ideal sheaf of the universal subscheme
$\Sigma_n\subset S^{[n]}\times S$, and $pr_{0i}$ is the projection
onto the product of the first factor $S^{[n]}$ and the $i$-th factor
of $S^l$. We shall denote $pr_0$ the projection onto the first
factor $S^{[n]}$, and $pr_i$  the projection onto the $i$-th factor
of $S^l$.

The induction statement, which will be proved by induction on $n$, is the following
generalization of Proposition \ref{prop2} (which is the $l=0$ case):
\begin{prop}
\label{induction}Let $ P\in CH(S^{[n]}\times S^l)$ be  a
 polynomial  expression  in $$
pr_0^*c_r(\mathcal{O}_{[n]}),\,pr_0^*c_s(T_n),\,pr_{0i}^*c_s(\mathcal{I}_n),\,1\leq i\leq l.$$
Then
for any $ \mu$ as above,
$$E_{\mu,l}^*(P)\in CH(S^m\times S^l)=CH(S^{m+l})$$
 is a polynomial
expression in the  $pr_j^*o,\,pr_{ik}^*\Delta,\,i,\,j,\,k\leq m+l$.
\end{prop}
{\bf Proof.}
Consider the smooth variety
$S^{[n,n-1]}$ parameterizing pairs  $(z,z')$ of subschemes of $S$, of length $n$ and $n-1$ respectively,
such that $z'\subset z$.

$S^{[n,n-1]}$ admits a natural map
$\rho$  to $S$, which to $(z,z')$ associates the residual point of
$z'$ in $z$. Together with the two natural projections $\psi$ to $S^{[n]}$ and
$\phi$ to $S^{[n-1]}$ respectively,
this gives two maps:
$$\psi:S^{[n,n-1]}\rightarrow S^{[n]},\,\,\sigma=(\phi,\rho):S^{[n,n-1]}\rightarrow S^{[n-1]}\times S.$$
$\sigma$ is birational; in fact it is the blow-up of
$S^{[n-1]}\times S$ along the incidence subscheme
$\Sigma_{n-1}\subset S^{[n-1]}\times S$. We shall denote by
$\mathcal{L}$ the line bundle $\mathcal{O}(-E)$ on $S^{[n,n-1]}$,
where $E$ is the exceptional divisor of $\sigma$. Thus we have
$$Im\,(\sigma^*\mathcal{I}_{n-1}\rightarrow  \mathcal{O}_{S^{[n,n-1]}})=\mathcal{O}(-E).$$
The map $\psi$ has degree $n$, and $(\psi,\rho)$ is a birational map
from $S^{[n,n-1]}$ to the incidence subscheme $\Sigma_{n}\subset
S^{[n]}\times S$.

Let now  $\mu=\{\mu_1,\ldots,\mu_m\}$ be a partition of
$\{1,\ldots,n\}$, and $S_\mu\cong S^m\subset S^n$ be as above.
Consider the fibered product
$$S_\mu\times_{S^{(n)}} S^{[n,n-1]},$$
which is also equal to
$$E_\mu\times_{S^{[n]}}S^{[n,n-1]}.$$
It obviously has exactly  $m$ components dominating $S_\mu$,
according to the choice of the residual point. Let us choose one
component, say the one where over the generic point
$(x_1,\ldots,x_n)\in S_\mu$, the residual point is $x_n$. Let $\mu'$
be the partition of $\{1,\ldots,n-1\}$ deduced from $\mu$ by putting
$$\mu'_i=\mu_i,\,\,{\rm if}\,\,n\not\in \mu_i,\,\,\mu'_i=\mu_i\setminus\{n\},\,\,{\rm if}\,\,n\in \mu_i.$$
Let us denote by $E_{\mu,\mu'}\subset S_\mu\times S^{[n,n-1]}$  the
underlying reduced variety of the component defined above, and note
that via the projection $\pi$ from $S_\mu$ to $S_{\mu'}$ (forgetting
the $n$-th factor), and the map $\sigma$, we get a natural map
$$\chi_{\mu'}=(\pi,\sigma)_{\mid E_{\mu,\mu'}}:E_{\mu,\mu'}\rightarrow E_{\mu'}\times S.$$
On the other hand, we have the natural map
$$\chi_{\mu}:=(Id_{S_\mu},\psi)_{\mid E_{\mu,\mu'}}:E_{\mu,\mu'}\rightarrow E_\mu.$$
Now, observe that  the following diagram is commutative:
$$\begin{matrix} E_\mu&\stackrel{\chi_\mu}{\leftarrow}&E_{\mu,\mu'}
&\stackrel{\chi_{\mu'}}{\rightarrow}& E_{\mu'}\times S\cr
p_{\mu}\downarrow&&&& \downarrow(p_{\mu'},Id)\cr
S_\mu&&\stackrel{\pi'}{\rightarrow}&& S_{\mu'}\times S
\end{matrix},
$$
where $p_\mu$ is the restriction to $E_\mu\subset S_\mu\times
S^{[n]}$ of the first projection, and similarly for $p_{\mu'}$, and
where $\pi':S_\mu\rightarrow S_{\mu'}\times S$ is given by
$(\pi,{pr_n}_{\mid S_\mu})$. Note also that both $\chi_\mu$ and
$\chi_{\mu'}$ are generically finite of degree $1$. Thus we have the
following equalities :
$$E_{\mu,\mu'}^*\circ \psi^*= E_\mu^*:CH(S^{[n]})\rightarrow CH(S_\mu),$$
$$\pi'_*\circ E_{\mu,\mu'}^*\circ \sigma^*=(E_{\mu'}\times \Delta_S)^*:CH(S^{[n-1]}\times S)
\rightarrow CH(S_{\mu'}\times S).$$
Similarly, for any integer $l$, we can consider the induced correspondence
$$E_{\mu,\mu',l}:=E_{\mu,\mu'}\times \Delta_{S^l}$$
between  $S_\mu\times S^l$ and $ S^{[n,n-1]}\times S^l$.
Then we have the formulas
\begin{eqnarray}
\label{compat5fev}E_{\mu,l}^*=E_{\mu,\mu',l}^*\circ (\psi,Id_l)^*
 :CH(S^{[n]}\times S^l)\rightarrow CH(S_\mu\times S^l),\\
 \label{compat25fev}
E_{\mu',l+1}^*=\pi'_{l*}\circ E_{\mu,\mu',l}^*\circ
(\sigma,Id_l)^*:CH(S^{[n-1]}\times S^{l+1}) \rightarrow
CH(S_{\mu'}\times S^{l+1}).
\end{eqnarray}
Here, $Id_l$ denotes the identity of $S^l$, and $\pi'_l$ is defined by
$$\pi'_l=(\pi',Id_l):S_\mu\times S^l\rightarrow
S_{\mu'}\times S^{l+1}.$$
 Furthermore, for any $\gamma\in CH(S^{[n,n-1]}\times S^{l})$,
 one has
 \begin{eqnarray}\label{nouveau27} \pi'_{l*}\circ E_{\mu,\mu',l}^*(
\gamma)=E_{\mu',l+1}^*((\sigma,Id_l)_*\gamma).
 \end{eqnarray}
 Indeed, this  follows from the fact that the correspondences
 $E_{\mu,\mu'}\subset S_\mu\times S^{[n,n-1]}$ and
 $E_{\mu'}\times\Delta_S\subset S_{\mu'}\times S\times  S^{[n-1]}\times S$
 satisfy the relation:
 \begin{eqnarray}\label{27oct}
 (\pi', Id_{S^{[n,n-1]}})_*(E_{\mu,\mu'})=(Id_{S_{\mu'}},\sigma,Id_S)^*
(E_{\mu'}\times\Delta_S)
\end{eqnarray}
in $CH(S_{\mu'}\times S\times S^{[n,n-1]})$
and similarly with $l>0$.
From (\ref{27oct}), we deduce that for
$\gamma\in CH(S^{[n,n-1]})$, one has
$$\pi'_{*}\circ E_{\mu,\mu'}^*(
\gamma)=(p_{S_{\mu'}\times S})_*((\pi', Id_{S^{[n,n-1]}})_*(E_{\mu,\mu'})\cdot p_{S^{[n,n-1]}}^*\gamma)$$
$$=(p_{S_{\mu'}\times S})_*((Id_{S_{\mu'}},\sigma,Id_S)^*
(E_{\mu'}\times\Delta_S)\cdot p_{S^{[n,n-1]}}^*\gamma)$$
$$=(p_{S_{\mu'}\times S})_*\circ (Id_{S_{\mu'}},\sigma,Id_S)_*((Id_{S_{\mu'}},\sigma,Id_S)^*
(E_{\mu'}\times\Delta_S)\cdot p_{S^{[n,n-1]}}^*\gamma)
$$
$$=(p_{S_{\mu'}\times S})_*((E_{\mu'}\times\Delta_S)\cdot
(Id_{S_{\mu'}},\sigma,Id_S)_*(p_{S^{[n,n-1]}}^*\gamma))$$
 $$=(p_{S_{\mu'}\times S})_*((E_{\mu'}\times\Delta_S)\cdot p_{S^{[n-1]}\times S}^*(\sigma_*\gamma))$$
 $$=
E_{\mu',1}^*(\sigma_*\gamma),$$
which proves (\ref{nouveau27}) for $l=0$. One argues similarly for $l>0$.

 From (\ref{nouveau27}), using the projection formula, one deduces
 that for any $$\alpha\in
CH(S^{[n,n-1]}\times S^l),\,\, \beta\in CH(S^{[n-1]}\times S^{l+1}),$$ one
has :
\begin{eqnarray}
\label{compatmenton} \pi'_{l*}\circ E_{\mu,\mu',l}^*(\alpha\cdot
(\sigma,Id_l)^*\beta)=E_{\mu',l+1}^*((\sigma,Id_l)_*\alpha\cdot\beta).
\end{eqnarray}

The key point is now the following formulas proved by Ellingsrud,
G\"ottsche, Lehn in \cite{EGL}: here we work on the $K_0$ groups
(the varieties considered are smooth and projective). The morphism
  $\phi^{!}:K_0(Y)\rightarrow K_0(X)$ for a morphism
$\phi:X\rightarrow Y$ between smooth varieties is induced by the morphism $\phi^*$ on vector bundles.
The morphism $M\mapsto M^{\vee}$ is induced by the morphism $E\mapsto E^*$ on vector bundles, and
the product $\cdot$ is induced by the tensor product between vector bundles.
Then we have (here we use for simplicity the fact that $K_S$ is trivial) :
\begin{theo}\label{theoEGL} (\cite{EGL}, Lemma 2.1 and Proposition 2.3)
We have in $K_0(S^{[n,n-1]})$ :
\begin{eqnarray}\label{T}
\psi^{!}T_n=\phi^{!}T_{n-1}+\mathcal{L}\cdot\sigma^{!}\mathcal{I}_{n-1}^{\vee}-\rho^{!}(1-T_S).
\end{eqnarray}
\begin{eqnarray} \label{On}
\psi^{!}\mathcal{O}_{[n]}=\phi^{!}\mathcal{O}_{[n-1]}+\mathcal{L}.
\end{eqnarray}
Furthermore, we have in $K_0(S^{[n,n-1]}\times S)$:
\begin{eqnarray} \label{In}(\psi,Id_S)^{!}\mathcal{I}_n=(\phi,Id_S)^{!}\mathcal{I}_{n-1}-
pr_0^{!}(\mathcal{L})\otimes( \rho,Id_S)^{!}\mathcal{O}_{\Delta_S}.
\end{eqnarray}
\end{theo}
Another very important property is
\begin{lemm} \label{li}(\cite{EGL}, Lemma 1.1)
In $CH(S^{[n-1]}\times S)$, we have the relation
$$ \sigma_*(c_1(\mathcal{L})^i)=(-1)^ic_i(-\mathcal{I}_{n-1}).$$
\end{lemm}
Theorem \ref{theoEGL} can be translated into statements concerning the Chern classes of the considered
sheaves (or elements of the $K_0$ groups).
Namely we conclude from (\ref{T}) that
the Chern classes $c_i(T_n)$ satisfy the property that
$\psi^*c_i(T_n)$ can be expressed as polynomials in
$$\phi^* c_j(T_{n-1}),\, c_1(\mathcal{L}),\,\sigma^*c_s(\mathcal{I}_{n-1}),\,\rho^*c_2(T_S)=24\rho^*o.$$
Similarly, we get from (\ref{On}) that
the Chern classes $c_i(\mathcal{O}_{[n]})$ satisfy the property that
$\psi^*c_i(\mathcal{O}_{[n]})$ can be expressed as polynomials in
$$\phi^* c_j(\mathcal{O}_{[n-1]} ),\, c_1(\mathcal{L}).$$
Finally, from (\ref{In}) we conclude that the Chern classes of
$\mathcal{I}_n$ satisfy the property that
$(\psi,Id_S)^*c_i(\mathcal{I}_{n})\in CH(S^{[n,n-1]}\times S)$ can
be expressed as polynomials in
$$(\phi,Id_S)^*c_j(\mathcal{I}_{n-1}),\,\,pr_0^*c_1(\mathcal{L}),\,\,(
\rho,Id_S)^*c_s(\mathcal{O}_{\Delta_S}) .$$ Note that because $K_S$
is trivial,  the Chern classes of $\mathcal{O}_{\Delta_S}$ reduce to
$$c_2(\mathcal{O}_{\Delta_S})=-\Delta_S \in CH^2(S\times S)$$ and
$c_4(\mathcal{O}_{\Delta_S})$, which is proportional to $(o,o)$ as
$c_2(T_S)$ is proportional to $o$.

Let now  $P \in CH(S^{[n]}\times S^l)$ be  a
 polynomial  expression  in $$
pr_0^*c_r(\mathcal{O}_{[n]}),\,pr_0^*c_s(T_n),\,pr_{0i}^*c_t(\mathcal{I}_n),\,1\leq
i\leq l$$ as in Proposition \ref{induction}. Applying
(\ref{compat5fev}), we get
\begin{eqnarray}\label{Pjeudi}E_{\mu,l}^*(P)=E_{\mu,\mu',l}^*\circ (\psi,Id_l)^*(P).
\end{eqnarray}
As just explained above, $(\psi,Id_l)^*(P)\in CH(S^{[n,n-1]}\times
S^l)$ can be expressed as a polynomial  in
$$(\phi,pr_{i})^*c_t(\mathcal{I}_{n-1}),\,pr_0^*c_1(\mathcal{L}), (\phi\circ
pr_0)^*c_r(\mathcal{O}_{[n-1]}),\,(\phi\circ pr_0)^*c_s(T_{n-1}),$$
$$(pr_{1,i}\circ(\sigma,Id_l))^*\Delta_S,\,(pr_i\circ(\sigma,Id_l))^*o,\,1\leq
i\leq l+1.
$$

 Observing that
$$\phi\circ pr_0:S^{[n,n-1]}\times
S^l \rightarrow S^{[n-1]}$$ is equal to $pr_0\circ(\sigma,Id_l)$,
the variables above can all be expressed as pull-back via
$(\sigma,Id_l)$ of the following variables in $CH(S^{[n-1]}\times
S^{l+1})$:
\begin{eqnarray}\label{123}pr_{1,i}^*\Delta_S,\,pr_i^*o,\,1\leq i\leq l+1,
\\ \nonumber
pr_{0i}^*c_t(\mathcal{I}_{n-1}),\, pr_0^*c_r(\mathcal{O}_{[n-1]}),\,
pr_0^*c_s(T_{n-1}), \end{eqnarray} except for
$pr_0^*c_1(\mathcal{L})$. Thus we have in $CH(S^{[n,n-1]}\times S^l)$:
\begin{eqnarray}\label{psilP}(\psi,Id_l)^*(P)=\sum_ipr_0^*c_1(\mathcal{L})^i(\sigma,Id_l)^*Q_i,
\end{eqnarray}
where $Q_i\in CH(S^{[n-1]}\times S^{l+1})$ is a polynomial
expression in the  variables (\ref{123}).

From (\ref{psilP}) and (\ref{Pjeudi}), applying
(\ref{compatmenton}), we deduce that
\begin{eqnarray}\label{derjeudimenton}\pi'_{l*}(E_{\mu,l}^*(P)=
\pi'_{l*}(E_{\mu,\mu',l}^*\circ (\psi,Id_l)^*(P))= \,\,\,\,\,\,\,\,\,\,\,\,\,\,\,\,\,\,\,\,\,\\ \nonumber
\pi'_{l*}(E_{\mu,\mu',l})^*(\sum_ipr_0^*c_1(\mathcal{L})^i(\sigma,Id_l)^*Q_i)
= E_{\mu',l+1}^*(\sum_iQ_i\cdot
(\sigma,Id_l)_*(pr_0^*c_1(\mathcal{L})^i)) .
\end{eqnarray}
Using Lemma \ref{li}, we find that
$(\sigma,Id_l)_*(pr_0^*c_1(\mathcal{L})^i))$ is a polynomial
expression in the $pr_{0j}^*c_s(\mathcal{I}_{n-1})$, and thus
$$\sum_iQ_i\cdot
(\sigma,Id_l)_*(pr_0^*c_1(\mathcal{L})^i)$$ is a polynomial
expression in the variables (\ref{123}). Applying induction on $n$
and the projection formula to the right hand side, we conclude that
$\pi'_{l*}(E_{\mu,l}^*(P))$ is a polynomial expression in the
variables
$$  pr_j^*o,\,pr_{ik}^*\Delta_S,\,i,\,j,\,k,\,\leq l+m .   $$

There are finally two cases to consider here, according to whether
$\mid\mu(n)\mid=1$ or $\mid\mu(n)\mid\geq2$, where $\mu(n)$ is the
element of the partition $\mu$ to which $n$ belongs (so $\mid
\mu(n)\mid$ is the multiplicity of $n$ in the diagonal $S_\mu$). In
the first case, we have
$$\pi':S_{\mu}\cong S_{\mu'}\times S,$$
 while in the second case, we have $\pi:S_\mu\cong S_{\mu'}$ and $\pi'$
 is the  embedding of
 $S_\mu \cong S^m$ in $S_{\mu'}\times S\cong S^{m+1}$ which is given by the diagonal on the last
 factor.
 In the first case, $\pi'$ being an isomorphism, we proved that
$E_{\mu,l}^*(P)$ is a polynomial in the variables
$pr_{ij}^*\Delta_{S},\,pr_k^*o$. In the second case, we get that
$pr_{\mu'}\circ \pi'$ is an isomorphism from $S_\mu$ to $S_{\mu'}$,
and applying $(pr_{\mu'},Id_l)$ to both sides of
(\ref{derjeudimenton}), we get the same conclusion.

This proves Proposition \ref{induction}, and thus also Proposition
\ref{prop2}. \cqfd
\begin{rema} It is presumably the case that Proposition \ref{prop2} could be obtained as a consequence
of the Bridgeland-King-Reid-Haiman equivalence of categories between
the derived category of $S^{[n]}$ and the derived category of
$\mathfrak{S}_n$-equivariant coherent sheaves on $S^n$ (see
\cite{bKR}, \cite{Haiman}), combined with results on equivariant
$K$-theory  of Vistoli \cite{vistoli}, and Riemann-Roch type
theorems by Toen \cite{toen}.

However, the explicit computation of the equivariant complex
associated to a given sheaf on $S^{[n]}$ is rather complicated. It
is done in \cite{scala} for $\mathcal{O}_{[n]}$, but not for $T_n$,
and the computation is  more difficult than the method of
\cite{EGL}, that we have been using here.
\end{rema}

\section{Case of the variety of lines of a cubic fourfold}

We shall use the following notations: the cubic fourfold will be denoted by $X$ and its Fano variety
of lines by $F$. $F$ is contained in the Grassmannian $G:=G(2,6)$ of lines in $\mathbb{P}^5$,
and we shall denote by
$$l\in CH^1(F,\mathbb{Z}),\,c\in CH^2(F,\mathbb{Z})$$
the  Chern classes of the  rank $2$ quotient bundle $\mathcal{E}$ induced on $F$.
Thus if
\begin{eqnarray}\label{Q}
\begin{matrix}P&\stackrel{q}{\rightarrow}& X\cr
p\downarrow&&\cr F&&
\end{matrix}
\end{eqnarray}
is the incidence diagram, $P$ is a $\mathbb{P}^1$-bundle over $F$,
and $\mathcal{E}=R^0p_*q^*(\mathcal{O}_X(1))$.

We shall denote by $H\in CH^1(X)$ the class $c_1(\mathcal{O}_X(1))$ and by
$h$ its pull-back to $P$, $h=q^*H$.

Let $I\subset F\times F$ be the incidence subvariety, which is  the codimension
$2$ subset of $F\times F$ defined as
\begin{eqnarray}I= (p,p)(q,q)^{-1}(\Delta_X),\label{I}
\end{eqnarray}
where $\Delta_X$ is the diagonal of $X$. Thus $I$ is the set of
pairs $(\delta,\delta')$ of intersecting lines. We shall denote by
the same letter $I$ the class of $I$ in $CH^2(F\times F)$.

We start the proof with a few remarks concerning the Chern classes
of $F$. As it is known that $F$ is symplectic holomorphic, one has
$T_F\cong \Omega_F$, and thus only the even Chern classes of $F$ can be non zero. We shall denote them by
$c_2,\,c_4$.
It is immediate  to compute that $c_2$ and $c_4$ can be written as polynomials in
$c$ and $l$.
Indeed $F\subset G$ is defined as the zero set of a section of the vector bundle
$S^3\mathcal{E}_G$ on $G$, and thus the normal bundle of $F$ in $G$ is isomorphic to
$S^3\mathcal{E}$. The normal bundle exact sequence then shows
that the Chern classes of $F$ are polynomials in $l,\,c$ and in the Chern classes
of $G$ restricted to $F$. But the later are also polynomials in $c$ and $l$, as are
 the restrictions  of  all cycles on
the Grassmannian.

Thus, in this case, Theorem \ref{maintheo}, 2) is equivalent to the following :
\begin{theo} \label{theocubic} Any polynomial expression in
$D\in CH^1(F)$ and $c\in CH^2(F)$ which vanishes in cohomology, vanishes in $CH(F)$.
\end{theo}

We observe first that there is no cohomological relation in degree
$4$ of the form above. Indeed, as $F$ is a deformation of a
$S^{[2]}$, one knows that
$$H^4(F,{\mathbb Q})\cong S^2H^2(F,{\mathbb Q}).$$
Thus there is only one cohomological relation of the form
$$[c_2(F)]=P,$$
where $P\in S^2H^2(F,{\mathbb Q})$. But this $P$ is non degenerate
because its kernel is a sub-Hodge structure of $H^2(F,{\mathbb
Q})^*$, which must be trivial because it is stable under deformation
of $F$, and in particular under a deformation for which $NS(F)$
becomes trivial. Thus there cannot be any relation of the form
$$[c_2(F)]=Q,$$
where $Q\in S^2(NS(F))$, because $NS(F) $ never generates
$H^2(F,{\mathbb Q})$.

Thus we only have to study relations in $H^6$ and $H^8$.
We first deal with the relations between $l$ and $c$ in degree $8$. There are obviously two such relations,
as $l^4,\,c^2,\,l^2c$ are all proportional in $H^8(F,\mathbb{Q})$.
Let us prove:
\begin{lemm}\label{lcrel}
 There exists a $0$-cycle $o\in CH^4(F)$, which is of degree $1$, and such that
$$ l^4,\,c^2,\,l^2c$$
are multiples of $o$ in $CH^4(F)$.
\end{lemm}

{\bf Proof.} We observe first that for generic $X$, there is one
surface $\Sigma$  of class $c$ which is a singular rational surface
(namely, its desingularization is rational). Indeed, surfaces in the
class $c$ are surfaces of lines of hyperplane sections of $X$. When
an hyperplane section $Y$ acquires a node $x$, its surface of lines
becomes birationally equivalent  to a symmetric product $S^2E_x$,
where $E_x$ is the curve of lines in $Y$ (or $X$) passing through
$x$ (see \cite{CG}). This curve of lines has genus $4$, and imposing
four ``independent'' supplementary nodes to $Y$ creates four nodes
on the curve $E_x$, which remains irreducible, so that the
normalization of $E_x$ becomes rational. In that case, the
desingularization of the surface of lines of $Y$ is rational. Now,
for generic $X$ it is easy to see that there exists such an
hyperplane section $Y$ with five independent nodes (which means that
the associated vanishing cycles are independent).

Of course, all  points of $\Sigma$ are rationally equivalent in $F$.
For some particular $X$, it might be that the surface $\Sigma$
degenerates to a non rational surface, but it still will remain true
that all the points of the degenerate surface $\Sigma$ are
rationally equivalent in $F$.

 We shall denote by
$o\in CH^4(F)$ this degree $1$ $0$-cycle. As  $c^2$ is supported on $\Sigma$,
$c^2$ is a  multiple of $o$ in $CH^4(F)$. Similarly $c\cdot l^2$ is supported on $\Sigma$,
hence it has  to be a  multiple of $o$ in $CH^4(F)$.

Next, with the same notations as above, we note that the curve $E_x$ is contained
in $\Sigma$. Thus we have a relation in $CH^4(X)$:
\begin{eqnarray} l\cdot E_x =\mu o,\label{1eqlundi}
\end{eqnarray}
for some coefficient $\mu$ equal to  the degree of $l\cdot E_x$.
The class of $E_x$ is computed as follows: As $ CH_0(X)=\mathbb{Z}$, this class does not depend
on $x$, and in fact we have :
$$3E_x=p_*h^4,$$
because $3x$ is rationally equivalent to $H^4$ in $X$.
Now we have the relation defining Chern classes:
$$(p^*l-h)h=p^*c$$
in $CH^2(P)$, which gives
$$ h^2=p^*l\cdot h-p^*c,\,\,h^3=p^*l\cdot(p^*l\cdot h-p^*c)-h\cdot p^*c=p^*(l^2-c)\cdot h-p^*(l\cdot c),$$
$$ h^4=p^*(l^2-c)\cdot(p^*l\cdot h-p^*c)-p^*(l\cdot c)\cdot h=
p^*(l^3-2lc)\cdot h-p^*((l^2-c)c).$$
Thus we have
\begin{eqnarray}\label{Ex}3E_x=l^3-2lc\,\,{\rm in}\,\, CH^3(F).
\end{eqnarray}
Equation (\ref{1eqlundi}) thus gives
a relation
$$ 3l(l^3-2lc)=\mu o,$$
and thus $l^4$ is also a multiple of $o$.

\cqfd
We now introduce a relation in the Chow ring of $F\times F$ which generalizes
the results obtained in \cite{voisinrelations} (which concerned the Chow ring of the surface of conics
of a Fano threefold). This relation will be essential to understand the group $CH_1(F)$.
\begin{prop} \label{propI}There is a quadratic relation in $CH^4(F\times F)$
\begin{eqnarray}\label{quadraticpropI}I^2=\alpha\Delta_F+\Gamma\cdot I +\Gamma',\end{eqnarray}
where $\alpha\not=0$, and $\Gamma$ is a codimension $2$ cycle of
$F\times F$ which is a degree $2$ polynomial in
$$l_1:=p_1^*l,\,l_2=p_2^*l,$$
and $\Gamma'$ is a codimension $4$ cycle which is a degree $2$
weighted polynomial in $l_1,\,l_2,\,p_1^*c,\,p_2^*c$.
\end{prop}
{\bf Proof.} We first prove the existence of a relation
of the above form, and  we will show later on
that the coefficient $\alpha$ is not $0$.

To get such a relation, it suffices to show the existence of a relation
\begin{eqnarray}\label{eqouvert}I^2_0=\Gamma\cdot I_0 +\Gamma'\,\,{\rm in}\,\,CH^4(F\times F\setminus\Delta_F),
\end{eqnarray}
where $\Gamma,\,\Gamma'$ are as above and $I_0$ is the restriction of $I$ to
$F\times F\setminus\Delta_F$.

Note that $I$ is the image in $F\times F$ via the map $(p,p)$ of
$$\widetilde{I}:=(q,q)^{-1}(\Delta_X).$$
Furthermore, over a point $(\delta,\delta')\in F\times F$, the fiber
of the map
$$p':=(p,p)_{\mid \widetilde{I}}:\widetilde{I}\rightarrow F\times F$$
identifies schematically to the intersection of the corresponding
lines $L,\, L'$ in $X$. Thus, away from the diagonal, this fibre is
a reduced point, and the restriction $p'_0$ of  $p'$ to
$\widetilde{I}_0:=\widetilde{I}\setminus (p')^{-1}(\Delta_F)$ is an
isomorphism onto $I_0$.

Furthermore, as $\widetilde{I}_0$ is a local complete intersection,
and $(p,p)$ is a submersion, $I_0$ is also a local complete
intersection, and thus $I^2_0$ is equal to $j_*(c_2(N_{I_0}))$,
where $N_{I_0}$ is the normal bundle of $I_0$ in $F\times
F\setminus\Delta_F$ and $j$ is the inclusion of $I_0$ in $F\times
F\setminus\Delta_F$. On the other hand, as $p'_0$ is an isomorphism
onto $I_0$, the normal bundle of $\widetilde{I}_0$ in $P\times P$
fits into a normal sequence
\begin{eqnarray}0\rightarrow {T_{P\times P/F\times F}}_{\mid \widetilde{I}_0}
\rightarrow N_{\widetilde{I}_0/P\times P}\rightarrow (p'_0)^*
N_{I_0/F\times F}\rightarrow0.
\end{eqnarray}

We deduce from this that ${p'_0}^* c_2(N_{I_0/F\times F})$ can be expressed as a polynomial in the
Chern classes $c_1,\,c_2$ of the normal bundle $N_{\widetilde{I}_0/P\times P}$
and in the Chern classes of ${T_{P\times P/F\times F}}_{\mid \widetilde{I}_0}$.

The later ones are polynomials in $h_1,\,l'_1,\,h_2,\,l'_2$, where
$$h_i=pr_i^*h,\,l'_i=pr_i^*(p^*l),\,i=1,\,2,$$
and $pr_i$ are the two projections of $P\times P$ onto $P$. Next we
observe that, as $\widetilde{I}=(q,q)^{-1}(\Delta_X)$, we have the
equalities
$$ c_i(N_{\widetilde{I}_0/P\times P})=q_0^*c_i(T_X),$$
where $q_0:\widetilde{I}_0\rightarrow X$ is the restriction of
$(q,q)$ to $\widetilde{I}_0$. But $c_i(T_X)$ are polynomials in $H$.
Thus we conclude that we have a relation:
$$I_0^2=p'_{0*}({p'}_0^* c_2(N_{I_0/F\times F}))$$
$$=p'_{0*}(P(h_i,\,l'_i)),$$
for some degree $2$ polynomial $P$ in $h_i,\,{l'_i}_{\mid
\widetilde{I}_0} $ (in fact $h_1=h_2$ on $\widetilde{I}_0$). This
can also be written as
$$I_0^2=(p,p)_*(P(h_i,l'_i)\cdot \widetilde{I})_{\mid F\times F\setminus\Delta_F}.$$

Let us now write the quadratic polynomial $P$ as
$$P=h_1A +h_2B+Q,$$
where $A,\,B$ are  linear in $h_i,l'_i$, and $Q$ is quadratic in $l'_1,\,l'_2$.
We have by the projection formula, noting that $l'_i=(p,p)^*l_i$,
$$(p,p)_*(Q(l'_i)\cdot \widetilde{I})=Q(l_i)\cdot I,$$
which is of the form $\Gamma'\cdot I$.

At this point we proved \begin{eqnarray} \label{ImentonF}
I_0^2=\Gamma'\cdot I_{0}+(p,p)_*((h_1A +h_2B)\cdot \widetilde{I})_{\mid
F\times F\setminus\Delta_F}.
\end{eqnarray} Finally, we observe that the diagonal of $X$ admits a
K\"unneth type decomposition:
$$\Delta_X=\Delta_1+\Delta_0,$$
where $\Delta_1$ can be written as a sum
$$\Delta_1=\sum_i\alpha_i H_1^i\cdot H_2^{4-i}$$
and $\Delta_0$ has the property that
\begin{eqnarray}\label{vanH}H_1\cdot \Delta_0=0,\,H_2\cdot \Delta_0=0\,\,{\rm in}\,\, CH^6(X\times X).
\end{eqnarray}
Here $H_i=pr_i^*H,\,i=1,\,2$, and $pr_i$ are the two projections on
$X\times X$. We obtain this decomposition as follows: we choose the
$\alpha_i$ in such a way that we have the following equalities
between intersection numbers:
$$\Delta_X\cdot H_1^i\cdot H_2^{4-i}=\Delta_1\cdot H_1^i\cdot H_2^{4-i},\,{\rm for}\,\, i=0,\ldots,4.$$
  Then
the cycle $\Delta_0=\Delta_X-\Delta_1$ is such that its image under
each inclusion
$$j_1:X\times X\hookrightarrow \mathbb{P}^5\times X,\,j_2:X\times X\hookrightarrow  X\times\mathbb{P}^5$$
is rationally equivalent to $0$, because
$j_{1*}\Delta_X={\Delta_{{\mathbb P}^5}}_{\mid {\mathbb P}^5\times
X}$. This implies (\ref{vanH}) because
$$j_1^*\circ j_{1*}=3H_1\cdot,\,\,\,\,j_2^*\circ j_{2*}=3H_2\cdot.$$

From the decomposition above, and recalling that
$$\widetilde{I}=(q,q)^{-1}(\Delta_X)=(q,q)^*\Delta_X,\,h_i=(q,q)^*H_i$$
we conclude that
$$h_1A\cdot \widetilde{I}=A\cdot(q,q)^*(H_1\cdot\Delta_X)=A\cdot(q,q)^*(H_1\Delta_1).$$
But as $H_1\Delta_1$ is a polynomial in $H_1,\,H_2$, it is then
clear that $(p,p)_*(h_1A\cdot\widetilde{I})$ is a cycle of the form
$\Gamma''$ as in the Proposition. Similarly for
$(p,p)_*(h_2B\cdot\widetilde{I})$. Thus, using (\ref{ImentonF}), the
existence of  a quadratic relation (\ref{quadraticpropI}) is proven.

 We now show that $\alpha\not=0$.
Mimicking the arguments in \cite{voisinrelations}, one sees that
there exist an hypersurface $W\subset F$ and a non zero coefficient
$\gamma\in \mathbb{Z}$ such that for each
 $\delta\in F$, there is a relation
$$\gamma \delta=S_\delta^2+z,$$
where $z$ is a $0$-cycle supported on $W$. Here $S_\delta$ is the
surface of lines of $X$ meeting $\delta$, so that
$S_\delta=I^*\delta$ in $CH^2(F)$ and
\begin{eqnarray}\label{Slcarre}S_\delta^2=\gamma \delta-z=(I^2)^*\delta\,\, {\rm in }\,\,CH^4(F).
\end{eqnarray}
 We have an equality
$$I^2=\alpha\Delta_F+ \Gamma\cdot I +\Gamma'\,\,{\rm in}\,\,CH^4(F\times F),$$
from which we deduce that $(I^2)^*$ acts as multiplication by
$\alpha$ on $H^{4,0}(F)\not=0$. On the other hand, (\ref{Slcarre})
together with the generalized Mumford theorem (cf
\cite{voisinbook}, Proposition 10.24 ), shows that $(I^2)^*$ acts as multiplication by
$\gamma$ on $H^{4,0}(F)$. Thus $\alpha=\gamma\not=0$.

\cqfd
We have the following corollary of Proposition \ref{propI}.
\begin{coro} \label{coro}Let $z\in CH_1(F)=CH^3(F)$ be a $1$-cycle.
Assume that $z$ is rationally equivalent to a combination of
rational curves $C_i\subset F$, $$z=\sum_in_iC_i,$$ that $z$ is
cohomologous to $0$, and that one (or equivalently any) point $x_i$
of  $C_i$ is rationally equivalent to  $o$ in $F$. Then $z=0$ in
$CH^3(F)$.
\end{coro}
{\bf Proof.} Indeed, observe that since $$z=\sum_in_iC_i,$$ with
$C_i$ rational, we have
\begin{eqnarray}\label{thu1}\Delta_{F*}z= \sum_i n_i(x_i\times
C_i+C_i\times x_i )\,\,{\rm in}\,\,CH(F\times F),
\end{eqnarray}
where $x_i$ is any point of $C_i$. Now $I^2$ is the restriction of
$I\times I$ to the diagonal $\Delta_{F\times F}$ of $F\times F$.
Thus we have $$(I^{2})^*z=({(I\times I)}^*(\Delta_{F*}z))_{\mid
\Delta_F}.$$
 From (\ref{thu1}), we conclude that
\begin{eqnarray}\label{thu2}
(I^{2})^*z=2\sum_in_iI^*C_i\cdot I^*x_i.
\end{eqnarray}
By assumption, we have $I^*x_i=I^*o$ in $CH^2(F)$, thus
(\ref{thu2}) is equal to
\begin{eqnarray}\label{thu3}
2I^*o\cdot \sum_in_iI^*C_i=2I^*o\cdot I^*z.
\end{eqnarray}
But $z$ is homologous to $0$, so $I^*z\in CH^1(F)$ is also homologous
to $0$, hence it is rationally equivalent to $0$. Thus
$(I^{2})^*z=0$ in $CH^3(F)$.

Now we apply Proposition \ref{propI} which gives a relation
$$\alpha z=(I^{2})^*z-(\Gamma\cdot I)^*z-\Gamma'^*z.$$
 As $(I^{2})^*z=0$, the right hand side   is equal to
$$-(\Gamma\cdot I)^*z-\Gamma'^*z.$$
But we know that both $I^*z$ and $l\cdot z$ are rationally
equivalent to $0$ : for the first, this was noticed just before, and
for the second, this is because it is a multiple of $o$ and
homologous to $0$. Hence it follows that $-(\Gamma\cdot
I)^*z-\Gamma'^*z=0$ and, as $\alpha\not=0$, we conclude that $z=0$.
\cqfd

As a consequence, we can start the computation of relations in
$CH^3(F)$ by
showing the following Lemma \ref{lethu1}:
Notice that $[l^3]$ and $[lc]$ are proportional in $H^6(F,{\mathbb Q})$.
Let this relation be $$[\mu cl-\nu l^3]=0\,\,{\rm in
}\,\,H^6(F,{\mathbb Q}),\,\mu\not=0,\nu\not=0.$$
\begin{lemm}\label{lethu1} We have the equality
$$  \mu cl-\nu l^3=0$$
in $CH^3(F)$.
\end{lemm}
{\bf Proof.} Indeed, it suffices to prove this relation for generic
$X$. In that case, we proved that the cycles $l^3$ and $lc$ are
supported on a rational surface of class $c$, all points of which
are rationally equivalent to $o$ in $F$. Thus the cycle $z=\mu
cl-\nu l^3$ satisfies the assumptions of Corollary \ref{coro}. \cqfd

In conclusion, we proved in Lemma  \ref{lcrel} and Lemma
\ref{lethu1} that all polynomial cohomological  relations between
$l$ and $c$ hold in $CH(F)$.

Let us decompose now $CH^1(F)$ as
$$CH^1(F)=<l>\oplus CH^1(F)_0,$$
where $CH^1(F)_0=p_*q^*CH^2(X)_{prim}$.
Recall the following from \cite{voisinbook}, 9.3.4.
Let $Z\in CH^2(X)_{prim}:=\{Z\in CH^2(F),\,[Z]\in H^4(X,{\mathbb
Q})_{prim}\}$. Write
$$q^*Z=hp^*D+p^*Z'.$$
Then from
$$H\cdot Z=0\,\,{\rm in}\,\,CH^3(X),$$
(see \cite{voisinbook}, 9.3.4), we get, using $h^2=hp^*l-p^*c$,
$$h^2p^*D+hp^*Z'=0=hp^*(lD+Z')-p^*(cD).$$
Thus we have  $D=p_*q^*Z$, and
\begin{eqnarray}\label{**}Z'=-lD,\,\,\,\,cD=0\,\,\,{\rm in}\,\,CH(F).
\end{eqnarray}
In particular
\begin{eqnarray}\label{eqthumidi}q^*Z=(h-p^*l)p^*D.
\end{eqnarray}

Let us deduce from this the following:
\begin{lemm} \label{lethumidiD2}For any $D\in CH^1(F)_0$, we have the
relations:
$$ l^2D^2=C q([D])o,$$
$$lD^2=C'q([D])E_x,$$
where $q$ is the Beauville-Bogomolov quadratic form on $H^2(F)$,
$C,\,C'$ are constants, and $E_x=p_*q^*x$ was already introduced and
shown to be proportional to  $l^3$ and $cl$ in $CH^3(F)$.
\end{lemm}
{\bf Proof.}
Note  that since $X$ is Fano, we have $CH^4(X)={\mathbb Q}$ and thus
\begin{eqnarray}\label{petite}Z^2=<Z,Z>x\end{eqnarray} for any $x\in X$.
Using (\ref{eqthumidi}), we get
\begin{eqnarray}\label{eqthumidi2}q^*(Z^2)=(h-p^*l)^2p^*(D^2).
\end{eqnarray}
Next we use the relations $cD=0$, $h^2=hp^*l-p^*c$, and (\ref{petite})  to rewrite
(\ref{eqthumidi2}) as
$$<Z,Z>q^*x=hp^*lp^*D^2-2hp^*lp^*D^2+p^*(l^2D^2)$$
$$=-hp^*(lD^2)+p^*(l^2D^2).$$
Note now that $<Z,Z>=-C'q([D])$ for some constant $C'$, as proved in
\cite{beaudonagi}, so that pushing forward via $p$ the above
expression, we get
$$C'q([D])E_x=lD^2.$$
Finally, applying $l$ to this, we get
$$l^2D^2=C'q(D)l\cdot E_x=Cq(D)o,$$
with $C=C'deg \,(l\cdot E_x)$. (We use (\ref{Ex}) and Lemma
\ref{lethu1} to get the last equality.) \cqfd Summing-up what we
have done up to now, we get:
\begin{prop}\label{thuprop}
Any polynomial relation $$[P]=0\,\,{\rm in}\,\,H^6(F,{\mathbb
Q})\,\,{\rm or\,\,in}\,\,H^8(F,{\mathbb Q}),$$
in the variables $l,\,c,\,D\in CH^1(F)_0$, which is of degree $\leq2$ in $D$,
is already satisfied in $CH^3(F)$, resp. $CH^4(F)$.
\end{prop}
{\bf Proof.} Indeed, consider first the case of $H^8$. The
polynomial expression $P$ is then of the form
$$P=cQ+l^2Q'+clA+l^3A'+\alpha c^2+\beta cl^2+\gamma l^4,$$
where $Q,\,Q'\in S^2 CH^1(F)_0$, $A,\,A'\in CH^1(F)_0$ and
$\alpha,\,\beta,\,\gamma$ are constants. But we know (cf (\ref{**}))
that
$$cQ=0,\,cA=0,$$
and that $l^2Q',\,c^2,\,cl^2,\,\gamma l^4$ are all multiples of $o$
(cf  Lemma \ref{lcrel}, Lemma \ref{lethumidiD2}).
On the other hand, as we proved that the cycle $l^3$ is rationally
equivalent to a cycle supported on a rational surface in the class $c$,
and all points of $\Sigma$ are rationally equivalent to $o$, it follows
that $l^3A'$ is also a multiple of $o$.
Thus $P$ is a multiple of $o$ in $CH^4(F)$, and as it is cohomologous to $0$,
it must be $0$.

Next we consider the case of degree $6$.
Then $P$ can be written as
$$P=lQ+cA+l^2A'+\alpha l^3+ \beta cl,$$
where $Q\in S^2CH^1(F)_0$, $A,\,A'\in CH^1(F)_0$ and $\alpha,\,\beta$
are constants.

We know that $cA=0$ and we proved already that the cycles
$$lQ,\,l^3,\,cl$$ are all proportional in $CH^3(F)$ (cf Lemma \ref{lethu1}, Lemma \ref{lethumidiD2}).
Using  these proportionality relations, we get an equality in
$CH(F)$:
$$P=l^2(A'+\gamma l),$$
where the number $\gamma$ depends on $Q,\,\alpha,\,\beta$
and involves the constants $\mu,\,\nu,\,C'$ of Lemma \ref{lethu1}, Lemma \ref{lethumidiD2}.
But we know that $[P]=0$, and thus the hard
Lefschetz theorem implies that $[A'+\gamma l]=0$. Thus, as   we are in
$CH^1(F)\subset H^2(F,\mathbb{Q})$, we have $A'+\gamma l=0$ and $P=0$.
\cqfd

We now turn to polynomials of degree at most $3$ in
$D$. Let us first consider the case of polynomials of degree $4$, that
is $P\in CH^4(F)$.
\begin{lemm}\label{frilem} Any polynomial expression
$P\in CH^4(F)$ in $l,\,c,\,D\in CH^1(F)_0$ which is of degree
at most $3$ in $D$ is a multiple of $o$. Thus, if
$[P]=0$ in $H^8(F,{\mathbb Q})$, then $P=0$.
\end{lemm}
{\bf Proof.} Indeed this was already proved for polynomial
expressions of degree at most $2$ in $D$ (cf Proposition
\ref{thuprop}), and thus, we only have  to consider expressions of
the form
$$P=lT,$$
where $T\in S^3CH^1(F)_0$. Now Lemma \ref{lethumidiD2} says  that for $D\in CH^1(F)_0$,
$lD^2 $ is proportional to $l^3$ in $CH^3(F)$. Hence $lD^3$ is
proportional to $l^3D$ in $CH^4(F)$. But by Proposition
\ref{thuprop}, we know that $l^3D$  is a multiple of $o$ in
$CH^4(F)$, as is any polynomial expression of degree $\leq2$ in $D$.\cqfd

We turn now to the cubic polynomial relations in $CH^3(F)$. First of
all we have the following lemma:
\begin{lemm}\label{cohD3} For any $D\in CH^1(F)_0)$, one has
\begin{eqnarray}\label{eqcohd3}
[D^3]= \frac{3}{q([l])}q([D])[l^2D].
\end{eqnarray}
\end{lemm}
{\bf Proof.} Recall from \cite{Ver}, \cite{Bo} that, in the complex cohomology algebra
$H^*(F,\mathbb{C})$, one has the relations
$$ {d'}^3=0,$$
for $d'\in H^2(F,\mathbb{C})$ such that $q(d')=0$.

It follows that we have more generally a relation of the form
$${d'}^3=q(d')A(d'),$$
where $A(d')\in H^6(F,\mathbb{C})$ is a linear function of $d'$.
We apply this to
$d'=d+\lambda [l]$, where $\lambda\in\mathbb{C}$, $d=[D],\,D\in CH^1(F)_0$.
Then we get, recalling that $q(d,[l])=0$,
\begin{eqnarray}\label{relmardi31}
d^3+3\lambda d^2[l]+3\lambda^2d[l]^2+\lambda^3[l]^3=(q(d)+\lambda^2q([l]))A(d')\,\,{\rm in}\,\,H^6(F,\mathbb{C}).
\end{eqnarray}
Write $A(d')=a(d)N+\lambda M$.
Then we get by taking the $0$-th order term in $\lambda$:
$$d^3=q(d)a(d)N.$$
 The order $2$ term in $\lambda$ gives now
 $$3d[l]^2=q([l])a(d)N,$$ from which
 we conclude that
 $$d^3=\frac{3}{q([l])}q(d)l^2d.$$

\cqfd
We will show the following proposition.
\begin{prop}\label{centerprop} For any $D\in CH^1(F)_0$, we have the relation
\begin{eqnarray}\label{centereqn}D^3=\frac{3}{q([l])}q([D])l^2D\,\,{\rm in}\,\,CH^3(F).
\end{eqnarray}
\end{prop}
Postponing the proof of Proposition \ref{centerprop},
we  conclude now the proof of Theorem \ref{theocubic}, or equivalently of Theorem
\ref{maintheo}, 2).

\vspace{0,5cm}

{\bf Proof of Theorem \ref{theocubic}.} Let us first treat the case
of a polynomial expression $P\in CH^3(F)$, which has to be  of
degree at most $3$ in $Pic\,F_0$. So assume $[P]=0$, where $P=T+lQ+
l^2 L+c L'+ C$, is the decomposition of $P$ into elements of
$Sym^\cdot CH^1(F)_0$ of degree $3$, $2$, $1$ and $0$ respectively,
whose coefficients are polynomials in $c,\,l$. We know from (\ref{**})
that $cL'=0$. We also know from Lemma \ref{lethumidiD2} and Lemma
\ref{lethu1} that $lQ$ and $C$ are proportional to $l^3$ in
$CH^3(F)$. Thus we have
$$lQ+C=\gamma l^3\,\,{\rm in}\,\,CH^3(F).$$
Finally, it follows from Proposition \ref{centerprop}  that $T$ is
equal in $CH^3(F) $ to   $l^2 D$ for some $D\in Pic\,F_0$.

Thus we have $P= l^2(D+L)+\gamma l^3$ in $CH^3(F)$ and the relation
$[P]=0$ implies
$$[l^2][D+L+\gamma l]=0 \,\,{\rm in}\,\,H^6(F,\mathbb{Q}).$$
But the hard Lefschetz theorem implies then that $[D+L+\gamma l]=0$.
 Thus $D+L+\gamma l=0$ and $P=0$.

To conclude the proof of the theorem, we now have  to consider the
case of a polynomial $P\in CH^4(F)$ of degree $4$ in $D\in
CH^1(F)_0$. But Proposition \ref{centerprop}  shows that, for any
$D\in CH^1(F)_0$, we have the relation
$$D^4=\frac{3}{q([l])}q([D])l^2D^2\,\,{\rm in}\,\,CH^4(F).$$
We proved in Lemma  \ref{lethumidiD2} that $l^2D^2$ is proportional
to $o$ in $CH^4(F)$. Thus $D^4$ is a multiple of $o$ and so is any
quartic homogeneous polynomial expression in $D\in Pic\,F_0$.

By Lemma \ref{frilem}, the same is true of any polynomial expression
of degree $\leq 3$ in $D$, with coefficients which are polynomials
in $l,\,c$. Thus any polynomial expression $P$ of degree $4$ in $D$,
with coefficients in $l,\,c$ is a multiple of $o$ in $CH^4(F)$. In
particular, if $[P]=0$, we have $P=0$.

\cqfd
{\bf Proof of Proposition \ref{centerprop}.} We first prove the result under the assumption that
$X$ contains no plane. We will show later on how to deduce the result  when $X$ contains planes.

Let us introduce the following object:
$$\widetilde{F}=\{(\delta_1,\delta_2)\in F\times F,\,\exists P\cong \mathbb{P}^2\subset\mathbb{P}^5,\,
P\cap X=2\delta_1+\delta_2\}.$$ Because we made the assumption that
$X$ does not contain any plane, $\widetilde{F}$ is irreducible, and
is the graph of the rational map $\phi: F\dashrightarrow F$
described in \cite{voisin}. We shall denote by
$$\tau:\widetilde{F}\rightarrow F,\,\tilde{\phi}:\widetilde{F}\rightarrow F,$$
the restrictions to $\widetilde{F}$ of the two projections. Thus
$\tau$ is birational and $\phi=\tilde{\phi}\circ\tau^{-1}$.

Note that $\widetilde{F}$ may be singular, which may imply that the
groups $CH_i(\widetilde{F})$ and $CH^{4-i}(\widetilde{F})$ differ,
and cause troubles because on one hand we  compute relations in
$CH_*(\widetilde{F})$, and on the other hand, we use intersection
product on $CH(\widetilde{F})$. However, there is a
desingularization of $\widetilde{F}$ which is obtained by a sequence
of blow-ups starting from $F$. We leave to the reader to adapt the
following arguments using this smooth model, and in the sequel, we
do as if $\widetilde{F}$ were smooth.

We will prove the following two Lemmas:
\begin{lemm}\label{lemmphiD} For $D\in CH^1(F)_0$, we have
$\tilde{\phi}^*D=-2\tau^*D$ in $CH^1(\widetilde{F})$.
\end{lemm}
\begin{lemm}\label{lemmphiI} Let $I\subset {F}\times F$ be the incidence subscheme
 defined in \ref{I}.
Then
\begin{eqnarray}\label{phiIdI}(\tilde{\phi},Id)^*I=-2(\tau,Id)^*I+Z
\end{eqnarray} in $CH^2(\widetilde{F}\times F)$, where $Z$ is a cycle of the
form \begin{eqnarray} \label{formulevalentin}Z=Z_1\times F+ D'\times
l +\widetilde{F}\times Z_2, \end{eqnarray} with $Z_1\subset
\widetilde{F}$ a codimension $2$ cycle, $D'\subset \widetilde{F}$ a
codimension $1$ cycle, $Z_2\subset F$ a
 codimension $2$ cycle.

\end{lemm}

Assuming these lemmas, let us show how to conclude the proof:
First of all, from Lemma \ref{lemmphiD}, we deduce that
for $D\in Pic\,F_0$, we have
\begin{eqnarray}\label{phid3}\tilde{\phi}^*D^3=-8\tau^*D^3\,\,{\rm in }\,\,CH^3(\widetilde{F}).
\end{eqnarray}

Next, from lemma \ref{lemmphiI}, we deduce that
\begin{eqnarray}\label{phiI2}(\tilde{\phi},Id)^*I^2=4(\tau,Id)^*I^2-4Z\cdot (\tau,Id)^*I+Z^2.
\end{eqnarray}
Note now that by definition of $\phi^*$ :
$$\tau_*\circ\tilde{\phi}^*=\phi^*,$$
acting on $CH(F)$.
Furthermore we have, applying $\tau_*,\,(\tau,Id)_*$ to (\ref{phid3}),
(\ref{phiI2}):
\begin{eqnarray}\label{phid3bis}\phi^*D^3=-8D^3,\\
\label{phiI2bis}(\phi,Id)^*I^2=4I^2-4I\cdot Z' +Z'',
\end{eqnarray}
where
$$Z':=(\tau,Id)_*Z,\,Z''=(\tau,Id)_*Z^2.$$

Observe now that
$$\phi^*((I^2)^*(z))=((\phi,Id)^*(I^2))^*(z),\,\forall z\in CH_1(F).$$
Combining this with (\ref{phiI2bis}) and the quadratic relation
(\ref{quadraticpropI}) given in Proposition \ref{propI}, we get, for
any $z\in CH_1(F)$: \begin{eqnarray}\label{forval2}\phi^*(\alpha
z+(\Gamma\cdot I)^*z +(\Gamma')^*z)=
 4(I^2)^*z -4(I\cdot Z')^*z+(Z'')^*z\\ \nonumber
=4(\alpha z+(\Gamma\cdot I)^*z +(\Gamma')^*z)-4(I\cdot
Z')^*z+(Z'')^*z. \end{eqnarray} Applying this to $z=D^3$ and  using
 (\ref{phid3bis}), we finally get
\begin{eqnarray}\label{relation27octobre}
-8\alpha D^3+\phi^*((\Gamma\cdot I)^*D^3 +(\Gamma')^*D^3)\\ \nonumber
=4\alpha D^3+4
((\Gamma\cdot I)^*D^3 +(\Gamma')^*D^3-(I\cdot
Z')^*D^3)+(Z'')^*D^3.
\end{eqnarray} In conclusion, we proved that
$$12\alpha D^3=\phi^*((\Gamma\cdot I)^*D^3 +(\Gamma')^*D^3)-
4((\Gamma\cdot I)^*D^3 +(\Gamma')^*D^3)-(I\cdot
Z')^*D^3)-(Z'')^*D^3.$$ We claim now that
$(\Gamma')^*D^3,\,\phi^*((\Gamma')^*D^3)$ and $(Z'')^*D^3$ are all
multiples of $l^3$ (or equivalently $cl$).

In the case of $(\Gamma')^*D^3$, this is a consequence of the fact
that $\Gamma'\in CH^4(F\times F)$ is a polynomial in
$pr_1^*l,\,pr_2^*l, \,pr_1^*c,\,pr_2^*c$, and of lemma \ref{lethu1}.
This implies also the claim for   $\phi^*((\Gamma')^*D^3)$, as one
shows easily (using Lemma \ref{lethu1}) that $\phi^*l^3$ is a
multiple of $l^3$. As for $(Z'')^*D^3$, we observe that we have for
any $z\in CH_1(F)$,
$$(Z'')^*z=\tau_*((Z^2)^*z),$$ and using formula (\ref{formulevalentin}) for
$Z$, this gives \begin{eqnarray}
\label{forval1}(Z'')^*z=2\tau_*(Z_1D') deg\,(l\cdot z). \end{eqnarray}
Hence it suffices to show that $\tau_*(Z_1D')$ is a multiple of
$l^3$. Now we have by (\ref{forval1}) applied to $l^3$:
$$(Z'')^*l^3=2\tau_*(Z_1D') deg\,l^4.$$
Thus it suffices to show that $(Z'')^*l^3$ is a multiple of $l^3$.
This follows now from (\ref{forval2}) applied to $z=l^3$, and from
the fact that
$$\phi^*l^3,\,(\Gamma\cdot I)^*l^3,\,(\Gamma')^*l^3,\, (I\cdot
Z')^*l^3$$ are all multiples of $l^3$. For the first three, this
follows easily from the definition of $\phi$ and from the form of
$\Gamma,\,\Gamma'$; for the last one, this follows from the fact
that, for any $z\in CH_1(F)$, $(I\cdot Z')^*z$ is a linear
combination of $\tau_*(Z_1)\cdot I^*(z)$ and $\tau_*D'\cdot
I^*(lz)$. Then the result is a consequence of the fact that
$$I^*l^3,\,I^*l^4,\,
\tau_*(Z_1),\,\tau_*D'$$
 are polynomial expressions in $l$ and $c$,
which is proved  using (\ref{phiIdI}) and the definitions of
$\widetilde{F}$ and $I$.

Next recall that the codimension $2$-cycle $\Gamma$ is a linear
combination of $l_1^2,\,l_2^2,\,l_1l_2$ on $F\times F$. Thus
$(\Gamma\cdot I)^*D^3$ is a combination of $l^2I^*(D^3)$ and of
$lI^*(lD^3)$. Next, for the same reason, $(I\cdot Z')^*D^3$ is a
linear combination of $\tau_*(Z_1)\cdot I^*(D^3)$ and $\tau_*D'\cdot
I^*(lD^3)$, that is of
$$l^2\cdot I^*(D^3),\,c\cdot I^*(D^3),\,l\cdot I^*(lD^3).$$

Thus our relation (\ref{relation27octobre}) becomes:
\begin{multline}
\label{dimancheder}12\alpha D^3=\phi^*(\mu l^2 I^*(D^3)+\nu l
I^*(lD^3))\\
 +\mu'l^2I^*(D^3)+\nu' lI^*(lD^3)+\mu''cI^*(D^3)+\mu'''
l^3.
\end{multline}

Recall from Lemma \ref{frilem} that $lD^3$ is proportional to $o$.
Thus $lI^*(lD^3)$ is proportional to $lI_o$ which is a multiple of
$l^3$ and $cl$ in $CH^3(F)$. Furthermore, we mentioned already that
$\phi^*(l^3)$ is also proportional to $l^3$.

Next we have
\begin{lemm} For some constant $\beta$, and for any $D\in CH^1(F)_0$, one has
$$ I^*(D^3)=\beta q([D])D.$$
\end{lemm}
{\bf Proof.} Indeed, as we are in $CH^1(F)$, it suffices to show this in $H^2(F,\mathbb{Q})$.
But we know that $[D^3]=\frac{3}{q(l)}q(D)[l^2D]$.
Thus it suffices to show that for some constant $\beta'$, and for any $[D]\in H^2(F,\mathbb{Q})_0$,
$$I^*([l]^2[D])]=\beta' [D]. $$
This is immediate because the left hand side is a morphism of Hodge structure
from $H^2( F,\mathbb{Q})_0$ to  $H^2( F,\mathbb{Q})$ which is defined for
general $X$, hence has to be a multiple of the identity, because the
Hodge structure on $H^2( F,\mathbb{Q})_0$  for general $X$ is simple with
$h^{2,0}=1$, while  $H^2( F,\mathbb{Q})=H^2( F,\mathbb{Q})_0+\mathbb{Q}[l]$.
\cqfd
From this lemma, we get in particular that $cI^*(D^3)=0$, and we
deduce from (\ref{dimancheder}) a relation:
$$12\alpha D^3=\phi^*(\mu q([D])l^2D)+ \mu'_1 q([D])l^2D +\nu' l^3.$$

 Furthermore, we recall that by Lemma \ref{lemmphiD}
$$\tilde{\phi}^*D=-2\tau^*D.$$
Hence it follows that
$$\phi^*(l^2D)=\tau_*(-2\tau^*D\tilde{\phi}^*l^2)=-2D\phi^*l^2.$$
It is easy  to verify that $\phi^* l^2$ is a combination of $l^2$
and $c$. As $cD=0$, we conclude that $\phi^*(l^2D)$ is a multiple of
$l^2D$. Thus we finally proved that we have a relation
\begin{eqnarray}\label{fineq}12\alpha D^3=\mu'' q([D])l^2D+\nu'l^3.
\end{eqnarray}
On the other hand, we know that we have the cohomological relation
$$[D^3]=\frac{3}{q(l)}q([D])[l^2D].$$
Using the hard Lefschetz theorem, and comparing with the cohomological relation
$$12\alpha [D^3]=\mu'' q([D])[l^2D]+\nu''[l^3]$$
deduced from (\ref{fineq}), we conclude that $\nu'=0$, and that
$$\frac{\mu''}{12\alpha}=\frac{3}{q(l)}.$$
This concludes the proof of Proposition \ref{centerprop} when $X$ contains no plane.

It remains to see how to do when $X$ contains a plane.
Let $D:=D_Z$ for some primitive class $[Z]\in H^4(X,\mathbb{Q})$.
In that case, either $[Z]$ is a multiple of  the primitive component
$[H]^2-3[\mathbb{P}]$  of the cohomology
class of a plane $\mathbb{P}\subset X$, or it is not. In the later case, one can show
by deformation theory that
a generic deformation of $X$ preserving the class $Z$ does not preserve any plane
contained in $X$. Then we know that (\ref{centereqn}) is satisfied
by $D_t\in Pic\,F_t$  for the generic member of a family of deformations
of the pair $(D,F)$. Thus it is also satisfied by $(D,F)$.

Thus it remains only to consider the case where $D=D_Z,\,[Z]=[H]^2-3[\mathbb{P}]$.
Thus $D=l-3D_\mathbb{P}$, where $D_\mathbb{P}$ is the divisor of
lines meeting $P$. But this case is easy because away from
the dual plane $\mathbb{P}^*\subset F$,
$D_\mathbb{P}$ is isomorphic via $p$ to
$\widetilde{D}:=q^{-1}(\mathbb{P})\subset P$. It follows that
the restriction
$(D_\mathbb{P})_{\mid D_\mathbb{P}}$ identifies  away from $\mathbb{P}^*$
as $det\,q^*N_{\mathbb{P}/X}-{T_{P/F}}_{\mid \widetilde{D}}$, that is
to the restriction of a combination of $h$ and $l$ to $\widetilde{D}$.
 From this, one deduces easily that
(\ref{centereqn}) is satisfied in $F\setminus\mathbb{P}^*$, and as
it is satisfied in cohomology, while
$$CH_1(\mathbb{P}^*)=H_2(\mathbb{P}^*,\mathbb{Z})=\mathbb{Z}\subset H^3(F,\mathbb{Q}),$$
it follows that it is satisfied as well on $F$.

Thus Proposition \ref{centerprop} is proved, modulo  Lemmas \ref{lemmphiD} and \ref{lemmphiI}.

\cqfd {\bf Proof of Lemma \ref{lemmphiD}.} Note that $\tau:
\widetilde{F}\rightarrow F$ is the contraction of a ruled divisor
$E$ to the surface $T$ of points $l\in F$ having the property that
there is a $\mathbb{P}^3_l\subset \mathbb{P}^5$ which is everywhere
tangent to the corresponding line $\Delta_l\subset X$. (One verifies
that $T$ is always a surface, and the fiber of $\tau$ over $l\in T$
identifies to the $\mathbb{P}^1$ parameterizing planes
$\mathbb{P}^2$ contained in $\mathbb{P}^3_l$ and containing
$\Delta_l$, because $X$ contains no plane.)

Thus for any divisor $D\in CH^1(F)$, there must be a relation
$$ \tilde{\phi}^*D=\tau^*D'+\sum_i\alpha_i E_i\,\,{\rm in}\,\,CH_3(\widetilde{F}),$$
where the $E_i$ are the irreducible components of $E$.
 Here the $\alpha_i$ are computed as
 $D\cdot \tilde{\phi}(E_{i,l})$, where $E_{i,l}$ is the fiber of
 $E_i$ over $l\in T_i$. (Here $T_i$ is the irreducible component of $T$ corresponding
 to $E_i$.)
 However, the curve $\tilde{\phi}(E_{i,l})$ is the family of lines contained in a cubic
 surface $S$ in $X$ which is singular along the line $\Delta_l$.
 Thus the surface in $X$ swept out by the lines parameterized by
$\tilde{\phi}(E_{i,l})$ is the cubic surface $S$, and for
$D=D_Z$, with $Z\subset X$ a cycle with primitive cohomology class, one has
 $$\alpha_i=-D\cdot \tilde{\phi}(E_{i,l})=<Z,S>=0.$$
 Thus we have
$$\tilde{\phi}^*D=\tau^*D'\,\,{\rm in}\,\,CH_3(\widetilde{F}),$$
and clearly $D'=\phi^*D\in CH^1(F)$. But the action of $\phi^*$ on
$CH^1(F)_0$ is the restriction of the action of $\phi^*$ on
$H^2(F,\mathbb{Q})_0:=p_*q^*H^4(X,\mathbb{Q})_{prim}$. This action
is multiplication by $-2$, because it is multiplication by $-2$ on
$H^{2,0}(F)$ (cf \cite{voisin}), and for general $X$ the Hodge
structure on $H^2(F,\mathbb{Q})_0$ is simple. Thus $D'=-2D$ and the
lemma is proven. \cqfd

{\bf Proof of Lemma \ref{lemmphiI}.} We
observe first that it suffices to prove the lemma for generic $F$,
because the family of $\tilde{F}$ parameterized by the set
$U\subset{\mathbb P}(H^0({\mathcal O}_{{\mathbb P}^5}(3)))$
corresponding to smooth cubic hypersurfaces which do not contain a
plane is flat.

Next we note that because $Pic^0F=0$, (which implies that divisors
on any product $K\times F$ are rationally equivalent to sum of
pull-backs of divisors on each factor,), and $Pic\,F={\mathbb Z}$,
which implies that divisors on $ F$ are rationally equivalent to a
multiple of $l$,  any codimension $2$ cycle in $\widetilde{F}\times
F$ which is supported on $D\times F$ is of the form
$$Z_1\times F+D'\times l,$$
where $Z_1$, $D'$ have respectively codimension $2$ and $1$ in
$\widetilde{F}$.

We use now the fact that for ${L}\in {F}$, the points ${L}$ and $
\phi({L})$ of $F$ parameterize lines
$$\Delta_{{L}},\, \Delta_{\phi({L})}$$
in  $X$ which satisfy the property
$$ 2\Delta_{L}+ \Delta_{\phi({L})}=H^3\,\,{\rm in}\,\,CH^3(X).$$
Thus we also have
$$ 2I_L+I_{\phi(L)}=C\,\,{\rm in}\,\,CH^2(F),$$
where $C=p_*q^*H^3$ is a constant. We then apply the Bloch-Srinivas
argument \cite{Blochsrinivas} (\cite{voisinbook},10.3.1), to conclude
that $2(\tau,Id)^*I+(\tilde{\phi},Id)^*I$ is rationally equivalent to the sum of  a
cycle of the form $\widetilde{F}\times C$ and of a
 cycle $W$ supported (via the first  projection) on a divisor  of $\widetilde{F}$.
 We can thus apply the remark above,
which gives
$$2(\tau,Id)^*I+(\tilde{\phi},Id)^*I=\widetilde{F}\times C+ Z_1\times F+D'\times l,$$
that is formula (\ref{formulevalentin}) with $Z_2=C$.

\cqfd


\begin{thebibliography}{99}
\bibitem{agv} D. Abramovich, T. Graber, A. Vistoli. Algebraic orbifold quantum products.
Orbifolds in mathematics and physics (Madison, WI, 2001), 1--24,
Contemp. Math., 310, Amer. Math. Soc., Providence, RI, 2002.
\bibitem{beauville} A. Beauville. Vari\'et\'es k\"ahleriennes dont la premi\`ere
classe de Chern est nulle, J. Differential Geometry 18 (1983),
755-782.
\bibitem{beauBB} A. Beauville. On the splitting of the Bloch-Beilinson filtration, in
{\it Algebraic cycles and motives} (vol. 2), London Math. Soc. Lecture Notes 344, 38--53;
Cambridge University Press (2007).
\bibitem{beaudonagi} A. Beauville, R. Donagi.
La vari\'{e}t\'{e} des droites d'une hypersurface
cubique de dimension 4, C.R. Acad. Sc. Paris 301, 703-706 (1985).
\bibitem{BeVo}A. Beauville, C. Voisin.
On the Chow ring of a K3 surface, J. Algebraic Geom. 13 (2004), no. 3, 417--426.
\bibitem{Blochsrinivas} S. Bloch, V. Srinivas.  Remarks on
 correspondences and algebraic cycles,
 Amer. J. of Math. 105 (1983) 1235-1253.
\bibitem{bo1} F. A. Bogomolov.  The decomposition of K\"{a}hler manifolds
 with a trivial canonical class. (Russian) Mat. Sb. (N.S.) 93(135) (1974), 573--575
\bibitem{bo2} F. A. Bogomolov.   Hamiltonian K\"{a}hlerian manifolds.
(Russian) Dokl. Akad. Nauk SSSR 243 (1978), no. 5, 1101--1104.
\bibitem{Bo} F. A. Bogomolov. On the cohomology ring of a simple
hyper-K\"{a}hler manifold (on the results of Verbitsky). Geom. Funct. Anal. 6 (1996), no. 4, 612--618.
\bibitem{bKR} T. Bridgeland, A. King, M. Reid.  Mukai implies McKay:
the McKay correspondence as an equivalence of derived categories,  J. Amer. Math. Soc. 14 (2001), no. 3, 535--554.
\bibitem{Cami}  M.A. de Cataldo, L. Migliorini: The Chow groups and the motive of the
 Hilbert scheme of points on a surface, Journal of algebra 251 (2002), 824-848.
 \bibitem{CG} H. Clemens, Ph. Griffiths. The intermediate jacobian of the cubic threefold,
  Ann. of Math. (2)  95  (1972), 281--356.
\bibitem{EGL} G. Ellingsrud, L. G\"{o}ttsche, M. Lehn. On the cobordism class of the Hilbert scheme of a surface.
Journal of Algebraic Geometry, 10 (2001), 81 - 100.
\bibitem{fango} B. Fantechi, L. G\"{o}ttsche. Orbifold cohomology for global quotients,
Duke Math J., vol. 117 (2003), 197-227.
\bibitem{Haiman} M. Haiman. Hilbert schemes, polygraphs, and the Macdonald positivity conjecture,
J. Amer. Math. Soc. 14 (2001), 941-1006.
\bibitem{LS} M. Lehn, C. Sorger. The Cup Product of the Hilbert Scheme for K3 Surfaces.
Inventiones mathematicae 152 (2003) 305 - 329.
\bibitem{scala} L. Scala. Ph. D. Thesis (2005).
\bibitem{Ver} M. Verbitsky. Cohomology of compact hyper-K\"{a}hler manifolds and
 its applications. Geom. Funct. Anal. 6 (1996), no. 4, 601--611.
 \bibitem{toen} B. Toen. Th\'{e}or\`{e}mes de Riemann-Roch pour les champs de Deligne-Mumford.
 (French)  $K$-Theory 18 (1999), no. 1, 33--76.
\bibitem{vistoli} G. Vezzosi, A. Vistoli. Higher algebraic $K$-theory
of group actions with finite stabilizers. Duke Math. J. 113 (2002), no. 1, 1--55.
\bibitem{voisin} C. Voisin.  Intrinsic pseudo-volume forms and K-correspondences,  Proceedings of the
 Fano Conference, (Eds A.Collino, A. Conte, M. Marchisio), Publication de l'Universit\'{e} de Turin (2004).
 \bibitem{voisinrelations} C. Voisin. Relations dans l'anneau de Chow de la surface des coniques des
vari\'et\'es de Fano,
 Amer. J. of Mathematics 112 (1990)
877-898.
\bibitem{voisinbook} C. Voisin. {\it Hodge Theory and Complex Algebraic Geometry}  II, Cambridge studies
in advanced mathematics 77, Cambridge Univ. Press (2003).
\end{thebibliography}
\end{document}